\documentclass[a4paper,leqno]{amsart}

\usepackage{amsmath, amssymb, amsfonts, euscript} 

\def\eps{\varepsilon}
\def\Csc{C_{\text{sc}}}
\def\Chf{C_{\text{hf}}}
\def\comp{_{\text{comp}}}
\def\loc{_{\text{loc}}}
\def\sD{\sqrt{-\Delta}}
\def\tOm{\theta_{\Omega}}
\def\tT{\theta_{T}}
\def\TOm{T_{\Omega}}

\def\vv{\underline{v}}
\def\ww{\underline{w}}
\def\ff{\underline{f}}
\def\kk{\underline{k}}
\def\aa{\underline{a}}
\def\rr{\underline{r}}

\def\Mt{\tilde{M}}
\def\nt{\tilde{n}}
\def\Deltat{\tilde{\Delta}}
\def\omj{{\omega_{j}}}
\def\ldk{\lambda_{k}}

\def\ldn{\lambda_{n}}

\def\sldn{\sqrt{\lambda_{n}}}

\def\otensprod{\otimes}
\def\ctensprod{\overline{\otimes}}
\def\Sset{\subset}

\DeclareMathOperator{\D}{D} 
\DeclareMathOperator{\supp}{supp}
\DeclareMathOperator{\dist}{dist}

\DeclareMathOperator{\id}{id}
\DeclareMathOperator{\Real}{Re}
\DeclareMathOperator{\Imaginary}{Im}
\renewcommand{\Re}{\Real}
\renewcommand{\Im}{\Imaginary}
\newcommand{\C}{\mathbb{C}}
\newcommand{\R}{\mathbb{R}}
\newcommand{\N}{\mathbb{N}}

\providecommand{\abs}[1]{\lvert #1 \rvert}
\providecommand{\norm}[1]{\lVert #1 \rVert}
\providecommand{\set}[1]{\left\{ #1 \right\}}
\newcommand{\Vect}[1]{\mathop{\rm Vect\left\{#1\right\}}}
\newcommand{\cd}[1]{\mathop{\bf #1}\nolimits}

\newcommand{\res}[1]{_{\rceil #1}}
\newcommand{\adh}[1]{\, \overline{#1}}
\newcommand{\conj}[1]{\overline{#1}}

 \theoremstyle{plain}
 \newtheorem{theorem}{Theorem}[section]
 \newtheorem{theorem*}{Theorem}
 \newtheorem{lemma}[theorem]{Lemma}
 \newtheorem{corollary}[theorem]{Corollary}
 \newtheorem{proposition}[theorem]{Proposition}

 \theoremstyle{definition}
 \newtheorem{definition}{{\sc Definition}}

 \theoremstyle{remark}
 \newtheorem{remark}[theorem]{{\sc Remark}}
 \newtheorem{remarks}[theorem]{{\sc Remarks}}

\begin{document}


\title[fast controls for Schr{\"o}dinger equation]{
How violent are fast controls \\
for Schr{\"o}dinger and plate vibrations~?}
\author[L. Miller]{Luc Miller}
\address{
{\'E}quipe Modal'X, JE 421, 
Universit{\'e} Paris X, B{\^a}t.~G,
200 Av.~de la R{\'e}publique,
92001 Nanterre, France.
}
\address{
Centre de Math{\'e}matiques, UMR CNRS 7640, 
{\'E}cole Polytechnique, 91128 Palaiseau, France.
}
\email{miller@math.polytechnique.fr}

\date{}

\thanks{This work was partially supported by the ACI grant
  ``{\'E}quation des ondes : oscillations, dispersion et contr{\^o}le''.
It was completed as the author was visiting 
the Mathematical Sciences Research Institute, Berkeley,
in the {\em Semiclassical Analysis} program. 
}


\begin{abstract}
Given a time $T>0$
and a region $\Omega$ on a compact Riemannian manifold $M$,
we consider the best constant, 
denoted $C_{T,\Omega}$, in the observation inequality
for the Schr{\"o}dinger evolution group 
of the Laplacian $\Delta$ with Dirichlet boundary condition:
$\displaystyle
\forall f\in L^{2}(M),\quad
\|f\|_{L^{2}(M)}
\leq C_{T,\Omega}
\|e^{it\Delta}f\|_{L^{2}((0,T)\times \Omega)}
$.
We investigate the influence of the geometry of $\Omega$
on the growth of $C_{T,\Omega}$ as $T$ tends to $0$.

By duality, $C_{T,\Omega}$  is also the controllability cost 
of the free Schr{\"o}dinger equation on $M$ with Dirichlet boundary condition 
in time $T$ by interior controls on $\Omega$.
It relates to hinged vibrating plates as well. 
We analyze the effects of wavelengths which are 
greater and lower than the control time $T$ separately.
We emphasize a tool of wider scope: {\em the control transmutation method}. 

We prove that $C_{T,\Omega}$ grows at least like $\exp(d^{2}/4T)$,
where $d$ is the largest distance 
of a point in $M$ from $\Omega$,
and at most like $\exp(\alpha_{*}L_{\Omega}^{2}/T)$,
where $L_{\Omega}$ is 
the length  of the longest generalized geodesic in $M$ 
which does not intersect $\Omega$,
and $\alpha_{*}\in]0,4[$ is the best constant 
in the following inequality for the Schr{\"o}dinger equation 
on the segment $[0,L]$ observed from the left end:
\linebreak 
$\displaystyle
\exists C>0, 
\forall f\in D(A), 
\|f\|_{H^{1}}
\leq C \exp(\alpha_{*} L^{2}/T)
\|\partial_{x} e^{itA}f{}_{\rceil x=L}
\|_{L^{2}(0,T)}
$,
where $A$ is the 
operator $\partial_{x}^{2}$ 
with domain $D(A)=\{f\in H^{2}(0,L) \,|\, Bf(0) = 0 = f(L) \}$
and the inequality holds with $B=1$ and with $B=\partial_{x}$.
We also deduce such upper bounds on product manifolds 
for some control regions which are not intersected by all geodesics.
\end{abstract}

\maketitle

\tableofcontents 


\section{Introduction}
\label{sec:intro}

\subsection{The problem}
\label{sec:pb}

Throughout the paper, 
$(M,g)$ is a smooth connected compact 
$n$-dimensional Riemannian manifold with metric $g$ 
and smooth boundary $\partial M$. 
When $\partial M\neq\emptyset$, $M$ denotes the interior 
and $\adh{M}=M\cup\partial{M}$.
Let $\dist:\adh{M}^{2}\to \mathbb{R}_{+}$ 
denote the distance function.
Let $\Delta$ denote the (negative) Dirichlet Laplacian on $L^2(M)$
with domain $H^1_0(M) \cap H^2(M)$.
Let $t\mapsto e^{-it\Delta}$ denote the 
Schr{\"o}dinger unitary group on $L^{2}(M)$.
The subset $\Omega$ of $M$ is always open.
For all positive time $T$, 
$\Omega_{T}$ denotes the time-space domain $]0,T[\times \Omega$.
Whenever generalized geodesics are mentionned (sections~4, 5 and 6), 
we make the additional assumptions that 
they can be uniquely continued 
at the boundary $\partial M$\footnote{\label{fnote:unique}
As in~\cite{BLR92}, to ensure this, we may assume either that 
$\partial M$ has no contacts of infinite order with its tangents
(e.g. $\partial M=\emptyset$),
or that $g$ and $\partial M$ are real analytic.
}.

\begin{definition}
\label{defin:cost}
For any $T>0$ and $\Omega\subset M$,
the {\em controllability cost} 
from $\Omega$ in time $T$
for the Schr{\"o}dinger equation on $M$ 
(with Dirichlet boundary condition if $\partial M\neq\emptyset$)
is the best constant, 
denoted $C_{T,\Omega}$, in the observation inequality:
\begin{equation}
\label{eq:cost}
\forall u_{0}\in L^{2}(M),\quad
\|u_{0}\|_{L^{2}(M)} 
\leq  C_{T,\Omega}
\|e^{-it\Delta}u_{0}\|_{L^{2}(\Omega_{T})} \ .
\end{equation}
\end{definition}

This observation inequality 
is a global and quantitative version of unique continuation 
from the domain  $\Omega_{T}$.
Let ${\bf 1}_{\Omega_{T}}$ 
denote the characteristic function of this space-time 
control region.
By duality (cf.~\cite{DR77}), 
the observation inequality $(\ref{eq:cost})$ is equivalent 
to the {\em exact controllability} of the free Schr{\"o}dinger equation
with Dirichlet boundary conditions
in time $T$ by interior controls on $\Omega$,
i.e. 
for all  
$u_{0}$ and $u_{T}$ in $L^{2}(M)$
there is a control function 
$g\in L^{2}(\mathbb{R}\times M)$
such that the solution $u\in C^{0}([0,\infty);L^{2}(M))$
(which can be defined by transposition) of:
\begin{equation} 
\label{eqSchr}
i\partial_{t}u - \Delta u=\cd{1}_{\Omega_{T}} g  
\quad {\rm in}\ \left]0,T\right[\times M, \quad 
u=0 \quad {\rm on}\ \left]0,T\right[\times\partial M,
\end{equation}
with Cauchy data
$u=u_{0}$ at $t=0$,
satisfies 
$u=u_{T}$ at $t=T$.
Moreover,  
$C_{T,\Omega}$ is also the best constant in the estimate:
\begin{gather*}
\|g\|_{L^{2}(\mathbb{R}\times M)}
\leq C_{T,\Omega}\|u_{0}-e^{iT\Delta}u_{T}\|_{L^{2}(M)}
\end{gather*}
for all data $u_{0}$, $u_{T}$, and all control $g$ 
solving this controllability problem.

This paper investigates the influence of 
the geometry of the control region $\Omega$ 
on the growth of the controllability cost $C_{T,\Omega}$
for the Schr{\"o}dinger equation 
as the control time $T$ tends to zero.
Fast controls of plate vibrations behave similarly since 
$\partial_{t}^{2} + \Delta^{2}=(\partial_{t} + i\Delta)(\partial_{t} - i\Delta)$
(precise statements for plates can be deduced straightforwardly
from our Schr{\"o}dinger results as in section~5 of~\cite{Leb92}).

\subsection{Main results}
\label{sec:res}

In subsection~\ref{sec:lb}, 
we deduce a finer statement of the following theorem 
(cf. theorem~\ref{th:lblf})
from a Gaussian estimate on the heat evolution for complex times 
(cf. proposition~\ref{prop:CGT}): 
\begin{theorem}
\label{th:lb}
The controllability cost of the Schr{\"o}dinger equation on $M$ 
from a nonempty subset $\Omega$ in short times 
(cf.\ definition~\ref{defin:cost}) 
satisfies the following geometric lower bound:
\begin{equation}
\label{eq:lb}
\liminf_{T\to 0} T \ln C_{T,\Omega} 
\geq \sup_{y\in M}\dist(y,\adh{\Omega})^{2}/4
\end{equation}
\end{theorem}

Our second result
concerns the most simple Schr{\"o}dinger controllability problem:
the Schr{\"o}dinger equation on a segment 
controlled at the right end through a Dirichlet condition.
Its generalization to Sturm-Liouville operators
(cf. theorem~\ref{th:1dSL}) 
is proved in section~\ref{sec:nonharm}
by the analysis of nonharmonic Fourier series.
This result is an upper bound of the same type as 
the lower bound in theorem~\ref{th:lb},
except that the rate $1/4$
is replaced by the technical rate
(resulting from lemma~\ref{lemM}):
\begin{equation}
\label{eqalphastar}
\alpha^{*}=4\left(\frac{36}{37}\right)^{2}
\ .
\end{equation}
The one dimensional boundary control version of 
the lower bound in theorem~\ref{th:lb} proves that 
the optimal rate $\alpha_{*}$ satisfies $\alpha_{*}\geq 1/4$. 
In its definition 
below, 
the notations for Sobolev spaces on the segment $[0,L]$ are: 
\begin{gather*}
H^{1}_{1}(0,L)=\{f\in H^{1}(0,L) \,|\, f(L)=0\}
 \mbox{ and }  
H^{1}_{0}(0,L)=\{f\in H^{1}_{1}(0,L) \,|\, f(0) = 0\}
\ .
\end{gather*}
\begin{definition}
\label{defin:alpha}
The rate $\alpha_{*}$ is the smallest positive constant 
such that the following controllability property holds: 
for all $\alpha > \alpha_{*}$
there exists $C>0$ such that, 
for all $k\in\{0,1\}$,
$L>0$, $T\in \left]0,\inf(\pi,L)^{2}\right]$
and $u_{0}\in H^{1}_{k}(0,L)$ 
the solution $u\in C^{0}([0,\infty);H^{1}_{k}(0,L))$
of the following Schr{\"o}dinger equation on $\left[0,L\right]$:
\begin{equation*} 
\label{eqSchr1d}
i\partial_{t}u - \partial_{s}^{2} u=0
\quad {\rm in}\ \left]0,T\right[\times \left]0,L\right[\, , \quad 
\partial_{s}^{k}u\res{s=0} = 0 = u\res{s=L}  \, ,\quad
u\res{t=0} = u_{0} \, ,
\end{equation*}
satisfies 
$\displaystyle
\|u_{0}\|_{H^{1}}
\leq C \exp(\alpha L^{2}/T)
\|\partial_{s} u_{\rceil s=L}
\|_{L^{2}(0,T)}
$.
\end{definition}
\begin{theorem}
\label{th:1d}
The rates in definition~\ref{defin:alpha} and $(\ref{eqalphastar})$ 
satisfy: $1/4\leq \alpha_{*}\leq \alpha^{*}< 4$.
\end{theorem}

Our third result, proved in section~\ref{sec:ub},
is an upper bound which is finite 
only under the {\em geodesics condition}\footnote{\label{fnote:geod}
In this context, 
this condition says that all 
generalized geodesics in $\adh{M}$ intersect the control region $\Omega$
(i.e. $L_{\Omega}<+\infty$ in theorem~\ref{th:ub}).
The {\em generalized geodesics} 
are continuous trajectories $t\mapsto x(t)$ in $\adh{M}$ 
which follow geodesic curves at unit speed in $M$
(so that on these intervals 
$t\mapsto \dot{x}(t)$ is continuous);
if they hit $\partial M$ transversely at time $t_{0}$,
then they reflect as light rays or billiard balls
(and $t\mapsto \dot{x}(t)$ 
is discontinuous at $t_{0}$); 
if they hit $\partial M$ tangentially 
then either there exists a geodesic in $M$
which continues $t\mapsto (x(t),\dot{x}(t))$
continuously 
and they branch onto it,
or there is no such geodesic curve in $M$
and 
then they glide at unit speed 
along the geodesic of $\partial M$
which continues $t\mapsto (x(t),\dot{x}(t))$
continuously until they may branch onto  a geodesic in $M$.
}
of C.~Bardos, G.~Lebeau and J.~Rauch, 
a.k.a. the geometric optics condition (cf.~\cite{BLR92}).
It is an application of the broader {\em control transmutation method}
(cf. sections~\ref{sec:backtransmut} and~\ref{sec:ub}). 
Here, it consists in writing the control $g$ for the Schr{\"o}dinger equation 
as a time integral operator applied to 
a control $f$ of the wave equation, 
i.e. $g(t,x)=\int_{\mathbb{R}} 
v(t,s)f(s,x)\, ds $,
where $f$ depends on $\Omega$ (not on $T$)
and the compactly supported kernel $v$ depends on $T$ and $L_{\Omega}$.
\begin{theorem}
\label{th:ub}
Let $\Omega\Sset M$ and let
$L_{\Omega}$ be the length of the longest generalized geodesic in $\adh{M}$
which does not intersect $\Omega$.
The controllability cost of the Schr{\"o}dinger equation 
from $\Omega$ in short times 
(cf.\ definition~\ref{defin:cost}) 
satisfies the following geometric upper bound 
(where $\alpha_{*}$ defined in definition~\ref{defin:alpha} 
satifies theorem~\ref{th:1d}):
\begin{equation}
\label{eq:ub}
\limsup_{T\to 0} T \ln C_{T,\Omega} \leq \alpha_{*}
L_{\Omega}^{2}
\ .
\end{equation}
\end{theorem}

Our last result is that 
the geodesics condition is not necessary 
for the controllability cost to grow at most like $\exp(C/T)$
as $T$ tends to $0$. 
In section~\ref{sec:ubprod},
a remark on the cost in an abstract tensor product setting
allows us to deduce from theorem~\ref{th:1d} and~\ref{th:ub} 
similar bounds in some settings violating the geodesics condition:
the boundary controllability of cylinders from one end  
(cf. theorem~\ref{th:cyl})
and the following semi-internal controllability on product manifolds 
(cf. theorem~\ref{th:prod} for the more abstract form).
\begin{theorem}
\label{th:ubprod}
Let $\Mt$ be a smooth complete $\nt$-dimensional Riemannian manifold
and $\Deltat$ denote the Laplacian on $L^2(\Mt)$
with domain $\{ u\in H^1_0(\Mt)\, |\, \Deltat u\in L^2(\Mt) \}$.
For all $T>0$ and all $\Omega \Sset M$,
the controllability cost $C_{T,\omega}$
of the Schr{\"o}dinger unitary group 
$t\mapsto e^{-it(\Delta +\Deltat)}$ on $L^{2}(M\times\Mt)$ 
from $\omega=\Omega\times \Mt$ in time $T$
is the controllability cost $C_{T,\Omega}$ 
of $t\mapsto e^{-it\Delta}$ on $L^{2}(M)$
from $\Omega$ in time $T$
(cf.\ definition~\ref{defin:cost}). 
In particular, with $\alpha_{*}$ and $L_{\Omega}$ as in theorem~\ref{th:ub}: 
$\limsup_{T\to 0} T \ln C_{T,\omega} \leq \alpha_{*}L_{\Omega}^{2}$.
\end{theorem}

\subsection{Background}
\label{sec:bg}

\subsubsection{Controllability for the Schr{\"o}dinger equation 
(and the plate equation)}

We survey from the geometric point of view 
the results on the exact controllability 
of the linear Schr{\"o}dinger equation in any positive time,
without discriminating 
boundary/interior observability/controllability 
for the Schr{\"o}dinger/plate equation. 
In this respect, the main result 
(proved by Lebeau in \cite{Leb92}) 
is  that the geodesics condition
is sufficient for boundary controllability 
on a smooth domain of $\mathbb{R}^{n}$ with a Riemannian metric
(cf. \cite{BZbb} for an alternative proof by resolvent estimates).
The same strategy applies 
to interior controllability 
(cf. the revisited proof in section~\ref{sec:hf})
and the control transmutation method yields yet another proof 
(cf. theorem~\ref{th:ub}). 

Further information on this condition
is obtained from the harmonic analysis of several examples
(some of them can be generalized and deduced directly from Lebeau's result, 
cf. section~\ref{sec:ubprod}).
The geodesics condition 
is not necessary for boundary controllability on a rectangle
(cf. \cite{KLS85}) 
and more generally on cylinders (cf. theorem~\ref{th:cyl}),
nor for interior controllability on a parallelepiped 
(cf. \cite{Har89}, \cite{Jaf90}, \cite{Kom92} in increasing generality),
on a torus (by the same proof), 
and more generally on a product manifolds (cf. theorem~\ref{th:prod}).
It is necessary for controllability on the sphere
(cf. \cite{Kom92})
except when the control region is an open hemisphere
(controllability holds in this case  
notwithstanding theorem~4.2 of \cite{Kom92}).

Burq also proved a controllability result 
(for a slightly more regular space of initial data) 
in the case of convex obstacles 
where the geodesics condition only fails for some 
hyperbolic trajectories of the geodesic flow.
Allibert studied the boundary control of revolution surfaces 
when the geodesics condition only fails for a single  
elliptic trajectory of the geodesic flow:
it can be checked that controllability in the natural spaces does not hold
(cf. section~2.1 in~\cite{All98}).
Recently,
Burq and Zworski proved in~\cite{BZbb} that controllability results 
for the classical and semiclassical Schr{\"o}dinger equation 
can be deduced from resolvent estimates, 
and give a striking application 
to the ergodic Bunimovich stadium. 

Other results assume geometric conditions 
which are more restrictive than 
the geodesics condition of Lebeau 
(and mostly stick to the Euclidean setting)
but require less smoothness than microlocal techniques:
they aim at more explicit estimates, nonlinear equations
and inverse problems.
The radial multiplier was used in 
\cite{Zua88}, \cite{Mac94}, \cite{Fab92} and~\cite{LT92}.
Carleman estimates are found 
in~\cite{Tat97}, \cite{TY99} (Riemannian setting), 
\cite{Zha01}, \cite{LTZ03}, \cite{BP02}.
Another approach based on local smoothing properties 
is sketched in \cite{LTay92} and \cite{HL96}. 

\subsubsection{Controllability cost}
The study of the controllability cost in short times 
was initiated by Seidman.
His first result in~\cite{Sei84} concerned the heat equation 
(see \cite{LMheatcost} for improvements and other references).
For many equations, the  controllability on a segment $[0,L]$ from one end
can be formulated as a window problem for series of complex exponentials
as in section~\ref{sec:1d}
(note that in this case $2L$ is the length of the longest generalized geodesic 
in $[0,L]$ which does not intersect one of the ends).
In~\cite{Sei86}, 
Seidman solved the window problem for purely imaginary exponentials 
corresponding to the Schr{\"o}dinger equation 
(he applied it to the plate equation in~\cite{KLS85})
therefore proving that in this setting the controllability cost grows 
at most like  $\exp(2(3\pi)^{2}\beta^{*}L^{2}/T)$ 
where $\beta^{*}\approx 4.17$
(or rather like $\exp(2\pi^{2}\beta^{*}L^{2}/T)$
if the sketchy remark~1 in section~4 works out).
Theorem~\ref{th:1d} improves on the constant appearing in this bound.
An example of Korevaar included in \cite{Sei86}
also proves that in this case the controllability cost grows
at least like $\exp(L^{2}/8T)$ 
(a computational slip in~\cite{Sei86} leads to $\exp(L^{2}/4T)$).
Seidman and his collaborators later 
treated the case of finite dimensional linear systems 
in~\cite{Sei88} and ~\cite{SY96},
and generalized the window problem to a larger class of complex exponentials 
in \cite{SG93} and \cite{SAI00}.

Phung's paper~\cite{Phu01} prompted our attention to the subject.
His theorem~2.3 proves that, under the geodesics condition, 
the cost of controlling data in $H^{1}_{0}(M)$
(one derivative more regular than in~theorem~\ref{th:ub})
grows at most as $\exp(C/T^{2})$
as $T$ tends to $0$ (one power of $T$ more than in~theorem~\ref{th:ub}
and no estimate on $C$).
Indeed, his one dimensional theorem~2.2
fell already short of the optimal power of $T$.
It can be checked that 
the usual Ingham theorem of harmonic analysis (for high frequencies) 
and the trick introduced by Haraux in~\cite{Har89}
(for the remaining low frequencies) 
yield the better (but still short of the optimal dependence in $T$) 
upper bound $\exp(C/(-T\ln T))$
(this is the approach followed in~\cite{JM01} for the wave equation).

\subsubsection{Transmutation}
\label{sec:backtransmut}
The strategy used by Phung to prove theorem~2.3 in~\cite{Phu01}, 
is what we have coined the {\em transmutation control method}.
Phung was inspired by~\cite{BdM74} and~\cite{KS96}
where the Schr{\"o}dinger semigroup on the whole space
is written as an integral over the wave group.
In fact, the method of transmutation applies between other kinds of equations 
(cf.~\cite{Her75} for a survey),
Kannai's formula being probability the best known example 
(cf. ~\cite{LMheatcost} for the corresponding application to heat control).

The most inspiring paper for both our lower and upper bound 
was~\cite{CGT82} which deduces 
geometric estimates on functions of the Laplace operators 
from the finite propagation speed of 
the even homogeneous wave group $W:s\mapsto \cos\left(s\sqrt{-\Delta}\right)$,
defined by: 
$w(s,x)=W(t)w_{0}(x)$ solves 
$\partial_{s}^{2}w - \Delta w=0$ 
in $\mathbb{R}\times M$ 
and $w=0$ on $\mathbb{R}\times\partial M$,
with Cauchy data
$(w,\partial_{s}w)=(w_{0},0)$ at $s=0$. 
It builds on the following transmutation formula 
which results from applying a spectral theorem 
to the Fourier inversion formula for an even function $F$:
\begin{equation}
\label{eqCGT}
F\left(\sqrt{-\Delta}\right)
=\int_{-\infty}^{+\infty}\hat{F}(s)W(s)\frac{ds}{2\pi}
\mbox{ , where } 
\hat{F}(s)=\int_{-\infty}^{+\infty}F(\sigma) \cos(s\sigma) d\sigma \ .
\end{equation}
When this formula is applied to evolution semigroups 
(like the Schr{\"o}dinger group $t\mapsto e^{it\Delta}$),
$F(\sigma)=\exp(tG(\sigma))$ where $t\geq 0$ is a time parameter
and $\hat{F}$ is a fundamental solution on the line
($\partial_{t}\hat{F}=\hat{G}{*}_{s}\hat{F}$ 
and $\hat{F}=\delta$ at $t=0$).
The {\em transmutation control method} consists in replacing 
this $\hat{F}$ 
by some  {\em fundamental controlled solution} on the segment $[-L,L]$ 
controlled at both ends.
We use the one dimensional theorem~\ref{th:1d} to 
construct this fundamental controlled solution
in subsection~\ref{sec:fcs}.
(Phung used a fundamental solution on the whole line 
controlled outside $[-L,L]$,
but it seems harder to estimate $\alpha_{*}$ for interior control.)
 
\subsection{The High/Low Frequencies issue}
\label{sec:plan}

Throughout the paper, 
$(\omega_{j})_{j\in \mathbb{N}^{*}}$
is a nondecreasing sequence of nonnegative real numbers 
and $(e_{j})_{j\in \mathbb{N}^{*}}$
is an orthonormal basis of $L^{2}(M)$
such that $e_{j}$ is an eigenvector of $-\Delta$
with eigenvalue $\omega_{j}^{2}$, i.e.:
\begin{equation}
\label{eqeigen}
-\Delta e_{j}=\omega_{j}^{2}e_{j} 
\quad\mbox{ and } \quad
e_{j}=0 \mbox{ on } \partial M
\ .
\end{equation}
The closed linear span of the vector set $\set{e_{j}}_{j\in J}$
is denoted by $\Vect{e_{j}}_{j\in J}$.

The spectral parameter $\omj$ can be considered as 
the frequency of the mode $e_{j}$.
For any given threshold $\mu>0$, 
the space of initial data can be decomposed into 
$L^{2}(M)=\Vect{e_{j}}_{\omj\leq\mu}\oplus\Vect{e_{j}}_{\omj>\mu}$
and this decomposition is invariant under the
Schr{\"o}dinger group $t\mapsto e^{-it\Delta}$.
The relevant notion of low (respectively high) frequencies 
in this paper correspond to wavelengths that are greater (respectively lower)
than the order of the control time, 
i.e. to $\mu\sim d/T$. 

Besides the main results already stated, 
the separate analysis of low and high frequencies
presented in sections~\ref{sec:lf} and~\ref{sec:hf}
give further insight into our initial problem. 
The cost of controlling low frequencies always grows like 
$\exp(C/T)$ as $T\to 0$. 
Under the geodesics condition $L_{\Omega}<+\infty$,
high frequencies are controlled at the much lower cost $C/\sqrt{T}$.
Though the upper bounds for low and high frequencies 
obtained respectively in subsections~\ref{sec:ublf} and~\ref{sec:ubhf} 
lead to conjecture the finiteness of 
$\limsup_{T\to 0} T \ln C_{T,\Omega}$
under the geodesics condition,
we emphasize that they do not suffice to prove this conjecture. 


Our high/low frequencies analysis leaves the following problem open:
can the controllability of the Schr{\"o}dinger equation
hold at a cost growing faster than $\exp(C/T)$
as $T$ tends to $0$~?
In others terms: are there $M$ and $\Omega\subset M$
such that $C_{T,\Omega}<+\infty$ for all $T>0$ 
and $\liminf_{T\to 0} T \ln C_{T,\Omega} = +\infty$~?
(n.b. theorem~\ref{th:ub} 
proves that violating the geodesics condition is necessary,
i.e. $\Omega$ must satisfy $L_{\Omega}=+\infty$.)
A positive answer would lead to the investigation 
of geometric conditions ensuring this ultra-violent behavior
(the examples of section~\ref{sec:ubprod} prove that 
violating the geodesics condition is not sufficient).


\section{Low frequencies}
\label{sec:lf}

In this section, we analyze how violent fast controls are  
for low frequency vibrations (cf. section~\ref{sec:plan}).

\subsection{Lower bound}
\label{sec:lb}

The purpose of this subsection is to prove the following refined version of 
theorem~\ref{th:lb}:
 
\begin{theorem}
\label{th:lblf} 
For all $\Omega\subset M$ and $d\in ]0,\sup_{y\in M}\dist(y,\adh{\Omega})[$: 
\begin{equation*}
\liminf_{T\to 0} T\ln 
\sup_{ u_{0}\in \dot{E}_{d/T} }
\frac{\|u_{0}\|_{L^{2}(M)}}{\|e^{it\Delta}u_{0}\|_{L^{2}(\Omega_{T})}}
\geq \frac{d^{2}}{4}
\mbox{, where } 
\dot{E}_{d/T}=\Vect{e_{j}}_{\omega_{j}\leq d/T }\setminus \{0\}
\ .
\end{equation*}
\end{theorem}

This lower bound follows from the construction of a very localized solution 
of the Schr{\"o}dinger equation with a large but finite number of modes.
For a short control time $T>0$, we consider a Dirac mass as far from $\Omega $
as possible, we smooth it out by applying the heat semigroup for a time $T/2$
and truncate frequencies larger than $d/T$,
and finally we take it as the data at time $T/2$ 
of the Schr{\"o}dinger equation.
The main ingredient in the proof of a similar bound for the heat equation 
(cf. theorem~2.1 in~\cite{LMheatcost})
was Varadhan's formula for the heat kernel in small time.
As a substitute here we prove:

\begin{proposition}
\label{prop:CGT} 
$\forall N\in\N$ such that $N>n/4$, $\forall \eps>0$, $\exists C_{\eps}>0$, 
$\forall d>0$, $\forall \Omega\subset M$, 
$\forall y\in M$ such that $\dist(y,\adh{\Omega})>d+\eps$,
$\forall z\in\C$ such that $\Re z >0$:
\begin{gather*}
\|e^{z\Delta}\delta_{y}\|_{L^{2}(\Omega)}\leq C_{\eps} 
\sqrt{\frac{|z|}{\Re z}}
\left( 1+(\Re z)^{-N}\left(1+\frac{d^{2}}{2|z|}\right)^{N}  \right)
\exp\left(-\frac{d^{2}\Re z}{4|z|^{2}}\right)
\ .
\end{gather*}
In particular, $\forall d>0$, $\exists B>0$, 
$\forall \Omega\subset M$, 
$\forall y\in M$ such that $\dist(y,\adh{\Omega})>d$,
$\forall z\in\C^{*}$, $|\Im z|\leq \Re z$ implies
$\|e^{z\Delta}\delta_{y}\|_{L^{2}(\Omega)}\leq 
B \exp\left(-d^{2}/(8\Re z)\right)$.
\end{proposition}

\begin{proof}
Our proof builds on the finite propagation speed 
and the boundedness on $L^{2}$ of 
the even homogeneous wave group $W:s\mapsto \cos\left(s\sqrt{-\Delta}\right)$
through the transmutation formula $(\ref{eqCGT})$.
Since $\Re z>0$ implies 
$e^{z\Delta}\delta_{y}\in 
\cap_{k\in\N}\D(\Delta^{k})\subset C^{\infty}(\adh{M})$
and $2N>n/2$ implies 
$\delta_{y}\in H^{-2N}_{\text{comp}}(M)\subset\D(\Delta^{N})'$
and therefore $W(s)\delta_{y}\in \D(\Delta^{N})' \subset \mathcal{D}'(M)$,
the following version of $(\ref{eqCGT})$ for $F(\sigma)=\exp(z\sigma^{2})$
makes sense for all $\varphi\in \D(\Delta^{N})$:
\begin{equation}
\label{eqCGTschr}
\left(e^{z\Delta}\delta_{y}, \varphi\right)_{L^{2}(M)}
=\int_{-\infty}^{+\infty}e^{-s^{2}/(4z)}f(s)\frac{ds}{\sqrt{4\pi z}}
\mbox{ , where } 
f(s)=\langle W(s)\delta_{y}, \bar{\varphi}\rangle
\ .
\end{equation}
As usual $\D(\Delta^{N})$ denotes the domain of the operator $\Delta^{N}$,
and $\D(\Delta^{N})'$ denotes it dual space
with respect to the duality product $\langle\cdot,\cdot\rangle$.
between distributions $\mathcal{D}'(M)$ and test functions 
$\mathcal{D}(M)=C^{\infty}_{\text{comp}}(M)$.

To prove the estimate in the proposition, 
we may assume $\varphi\in \mathcal{D}(\Omega)$
since  $\mathcal{D}(\Omega)$ is dense in $L^{2}(\Omega)$.
We deduce $f\in C^{\infty}(\R)$ with 
$\supp f \subset \left\{s\in\R | \, |s|\geq \dist(y,\adh{\Omega})\right\}$
since $\partial_{s}^{2}W\delta_{y}=\Delta W\partial(y)$
and $\supp W(s)\delta_{y}\subset  
\left\{x\in M | \, \dist(x,y)\leq |s|\right\}$
(propagation at unit speed).
Moreover 
$f(s)=\langle (1-\Delta)^{N}W(s)(1-\Delta)^{-N}\delta_{y},\varphi\rangle
=(1-\partial_{s}^{2})^{N}g(s)$
where $g(s)=\langle W(s)(1-\Delta)^{-N}\delta_{y}, \bar{\varphi}\rangle$
satisfies $\|g\|_{L^{\infty}}
\leq \|(1-\Delta)^{-N}\delta_{y}\|_{L^{2}}\|\varphi\|_{L^{2}}$
since $W(s)$ is bounded on $L^{2}$.
Therefore, we may integrate by parts and obtain:
\begin{equation}
\label{eqCGT1}
\begin{split}
\int e^{-s^{2}/(4z)}f(s)ds 
& = \int e^{-s^{2}/(4z)}(1-\partial_{s}^{2})^{N}g(s)ds \\
& = \int g(s) (1-\partial_{s}^{2})^{N}
\left(\chi_{\eps}(s)e^{-s^{2}/(4z)}\right) ds 
\end{split}
\end{equation}
where $\chi_{\eps}(s)=\chi(\frac{|s|-d}{\eps})$ 
is a smooth non-negative cut-off function satisfying $\chi_{\eps}(s)=1$
on $|s|\geq d+\eps$ and $\chi_{\eps}(s)=0$ on $|s|\leq d$.
From the simple estimate:
$\forall k\in\N$, $\exists C_{k}>0$,
\begin{gather*}
\forall z\in \C, \Re z>0:\ 
\left|\partial_{s}^{2k}e^{-s^{2}/(4z)}\right|
\leq  
\frac{C_{k}}{|z|^{k}}\left(1+\frac{s^{2}}{4|z|}\right)^{k}
\exp\left(-\frac{s^{2}\Re z}{4|z|^{2}}\right)
\ ,
\end{gather*}
we deduce, setting $\tau=d\sqrt{\Re z}/(2|z|)$:
\begin{equation}
\label{eqCGT2}
\begin{split}
\left|
\int\limits_{d}^{+\infty} \partial_{s}^{2k}e^{-s^{2}/(4z)} ds
\right|
& \leq \frac{C_{k}}{|z|^{k}}
\int\limits_{\tau}^{+\infty}
\left(1+s^{2}\frac{|z|}{\Re z}\right)^{k} e^{-s^{2}} \frac{d}{\tau}ds \\
& = \frac{C_{k}de^{-\tau^{2}}}{|z|^{k}\tau}
\int\limits_{0}^{+\infty}
\left(1+(s+\tau)^{2}\frac{|z|}{\Re z}\right)^{k} e^{-s^{2}}e^{-2s\tau}ds \\
& \leq  \frac{C_{k}de^{-\tau^{2}}}{|z|^{k}\tau}
\left(1+\frac{2\tau^{2}|z|}{\Re z}\right)^{k}\frac{|z|^{k}}{(\Re z)^{k}}
\int\limits_{0}^{+\infty}(1+2s^{2})^{2k}e^{-s^{2}}ds \\
& =\frac{C_{k}'|z|}{\sqrt{\Re z}}
(\Re z)^{-k}\left(1+\frac{d^{2}}{2|z|}\right)^{k}  
\exp\left(-\frac{d^{2}\Re z}{4|z|^{2}}\right)
\ .
\end{split}
\end{equation}

Equations $(\ref{eqCGTschr})$, $(\ref{eqCGT1})$ and $(\ref{eqCGT2})$
imply the first estimate in proposition~\ref{prop:CGT} with a $C_{\eps}$
which only depends on $C_{k}'$ and on $\eps$ 
through $\sup_{k\leq N}\|\partial^{k}\chi_{\eps}\|_{L^{\infty}}$.
Since $|\Im z|\leq \Re z$ implies $|z|^{2}\leq 2(\Re z)^{2}$,
the second estimate results from the first 
(at the expense of the uniformity in $d$).
\end{proof}

\begin{proof}[Proof of theorem~\ref{th:lblf}]
We shall use Weyl's asymptotics for eigenvalues:
\begin{equation}
  \label{eqWeyl}
  \exists W>0,\  
\#\{j\in \mathbb{N}^{*}\, | \, \omega_{j}\leq \omega \} \leq W \omega^{n}  
\end{equation}
and the following consequence of Sobolev's embedding theorem:
\begin{equation}
  \label{eqSobolev}
  \exists E>0,\,  \forall  j\in \mathbb{N}^{*},\
\|e_{j}\|_{L^{\infty}}\leq  E \omega_{j}^{n/2}
\end{equation}
(cf. section 17.5 in \cite{HIII} for example).
The unique continuation property for elliptic operators implies that 
$Y=\{ y\in M\setminus \adh{\Omega} \, |\, e_{1}(y)\neq 0 \}$
is an open dense set in $M\setminus\adh{\Omega}$,
so that the supremun in theorem~\ref{th:lblf}
can be taken over $y\in Y$ instead of $y\in M$.

Let $y\in Y$ and $D<d<\dist(y,\adh{\Omega})$
be fixed from now on.
Applying proposition~\ref{prop:CGT} with $z=T-it$ 
yields 
a positive constant $B$ such that: 
$\forall T>0$, 
$\forall t\in ]0,T]$, 
\begin{equation}
  \label{eqk}
\|e^{(T-it)\Delta}\delta_{y}\|_{L^{2}(\Omega)}
\leq B e^{-d^{2}/(8T)}
\ .
\end{equation}

Therefore, for all $T>0$,
we take as initial data the following finite modes approximation of 
$e^{T\Delta}\delta_{y}$:
$u_{0}^{T}(x)=\sum_{2T\omega_{j}\leq d}
\exp(-T\omega_{j}^{2})e_{j}(y)e_{j}(x)$,
and we are left with comparing $e^{(T-it)\Delta}\delta_{y}$
to the corresponding solution 
\begin{gather*}
u^{T}(t,x)=\left( e^{-it\Delta}u_{0}^{T} \right)(x)
=\sum_{2T\omega_{j}\leq d}
\exp((it-T)\omega_{j}^{2})e_{j}(y)e_{j}(x)
\ .
\end{gather*}
Using the unitarity of the Schr{\"o}dinger group on $L^{2}(M)$,
Parseval's identity and (\ref{eqSobolev}), 
we obtain
\begin{multline*}
\sup_{t\in ]0,T]} \| e^{(T-it)\Delta}\delta_{y}-u^{T}(t,x)\|_{L^{2}(M)}
\leq  \| e^{(T-iT)\Delta}\delta_{y} -u^{T}_{0}(x)\|_{L^{2}(M)} \\
=\sum_{2T\omega_{j}>d }
|e^{-T\omega_{j}^{2}}e_{j}(y)|^{2}
\leq E \sum_{2T\omega_{j}\geq d}
e^{-d\omega_{j}/2}\omega_{j}^{n}
\leq 
E'\sum_{2T\omega_{j}\geq d}
e^{-D\omega_{j}/2}
\ ,
\end{multline*}
for some $E'>0$.
But, Weyl's law (\ref{eqWeyl}) yields, 
for $c\geq c_{0}>0$ and $\gamma\geq \gamma_{0}>0$,
\begin{multline*}
\sum_{\omega_{j}\geq c}e^{-\gamma \omega_{j}}
=\sum_{k\in \mathbb{N}^{*}} 
\sum_{kc\leq \omega_{j} <(k+1)c}e^{-\gamma \omega_{j}}
\leq W \sum_{k\in \mathbb{N}^{*}} 
\left( (k+1)c \right)^{n} e^{-kc\gamma } \\
\leq  W_{\gamma_{0}} 
\sum_{k\in \mathbb{N}^{*}} e^{-kc\gamma}e^{(k+1)c\gamma/4}
= W_{\gamma_{0}} e^{-c\gamma/2}
\sum_{k\in \mathbb{N}} e^{-3kc\gamma/4} 
\leq W_{c_{0},\gamma_{0}} e^{-c\gamma/2}
\end{multline*}
where $W_{\gamma_{0}}$ and $W_{c_{0},\gamma_{0}}$
are positive real numbers which depend on their indexes 
but not on $c$ and $\gamma$.
Hence, with $c=d(2T)^{-1}>d/2=c_{0}$ and $\gamma=\gamma_{0}=D/2$, we obtain: 
\begin{equation*}
 \exists B'>0,\ 
 \forall t\in ]0,T]\ 
\| e^{(T-it)\Delta}\delta_{y}-u^{T}(t,x)\|_{L^{2}(M)}
\leq B' e^{-dD/(8T)}
\end{equation*}
Together with $(\ref{eqk})$,
this estimate yields,
setting $B''=B + B'$, for all $T>0$:
\begin{gather*}
\| u^{T} \|_{L^{2}(\Omega_{T})}
\leq \sqrt{T}Be^{-d^{2}/(8T)}
+ \sqrt{T}B'e^{-dD/(8T)}
\leq \sqrt{T}B''e^{-D^{2}/(8T)}
\ .
\end{gather*}
But using Parseval's identity and $y\in Y$, we have
for all $T\in \left] 0,1 \right]$:
\begin{gather*}
\|u_{0}^{T}\|_{L^{2}(M)}
=\left( 
\sum_{2T\omega_{j}\leq d}
|e^{-T\omega_{j}^{2}}e_{j}(y)|^{2}
\right)^{1/2}
\geq e^{-\omega_{1}^{2}} |e_{1}(y)|
>0
\ .
\end{gather*}
Hence, 
with $A=e^{-\omega_{1}^{2}} |e_{1}(y)|B''$ independent of $T$,
$u_{0}^{T}\in \dot{E}_{d/(2T)}$ satisfies 
\begin{gather*}
\forall T\in ]0, 1],\ 
\|u^{T}\|_{L^{2}(\Omega_{T})}
\leq A\sqrt{T}e^{-D^{2}/(8T)}\|u_{0}^{T}\|_{L^{2}(M)}
\ .
\end{gather*}
Applying the last estimate forward and backward in time yields
an estimate on a time interval of length $2T$:
$\|u^{T}\|_{L^{2}(]-T,T[\times \Omega)}\leq 
A\sqrt{2T}e^{-D^{2}/(4(2T))}\|u_{0}^{T}\|_{L^{2}(M)}$.
Since $\|u^{T}\res{t=-T}\|_{L^{2}(M)}=\|u_{0}^{T}\|_{L^{2}(M)}$
and $D<d$ is arbitrary,  
changing $2T$ into $T$ ends the proof of theorem~\ref{th:lblf}.
\end{proof}

\begin{remarks}
In the case $M=S^{1}$ (the unit circle),
the transmutation formula $(\ref{eqCGT})$ 
is essentially the Poisson summation formula.
In this sense, our construction is an extension of 
Korevaar's one dimensional example in~\cite{Sei86}.

Following~\cite{CGT82}, 
we could also prove point wise Gaussian estimates 
of the Heat kernel for complex times. 
Proposition~\ref{prop:CGT} is a short path to the estimate 
required by our construction.

For $z=h+ith$, this proposition is an analogue on the compact manifold $M$
of the localization estimate satisfied by the solution of the 
semiclassical  Schr{\"o}dinger equation $ih\partial_{t}u-h^{2}\Delta u=0$
in $\mathbb{R}^{n}$ with initial data $u_{0}(x)=\exp(-(x-y)^{2}/(4h))$,
i.e. a semiclassical coherent state centered at $y$ with no momentum.
\end{remarks}

\subsection{Upper bound at low frequencies}
\label{sec:ublf}

Carleman estimates are the most versatile tool to control low frequencies 
as epitomized by their application in theorem~3 in~\cite{LZ98} 
(also theorem~14.6 in~\cite{JL99}):

\begin{theorem}[\cite{LZ98},\cite{JL99}]
For all non-empty open subset $\Omega$ of $M$:
\begin{gather*}
\exists C>0,  \forall v\in \C^{\N^{*}}, \forall \mu>0,
\quad
\sum_{\omega_{j}\leq \mu} |v_{j}|^{2}
\leq Ce^{C\mu} \int_{\Omega} 
\left| \sum_{\omega_{j}\leq \mu} v_{j}e_{j}(x) \right|^{2}dx
\ .
\end{gather*} 
\end{theorem}

Applying this observation inequality for fixed time 
and integrating on $[0,T]$ yields:
\begin{gather*}
\exists C, \forall \mu>0, 
\forall u_{0}\in \Vect{e_{j}}_{\omega_{j}\leq \mu},
\quad
\|u_{0}\|_{L^{2}(M)}\leq 
\frac{C}{\sqrt{T}}e^{C\mu}\|e^{it\Delta}u_{0}\|_{L^{2}(\Omega_{T})}
\ .
\end{gather*}
As a counterpart to theorem~\ref{th:lblf}, 
taking $\mu=d/T$, we state:
\begin{corollary}
\label{cor:ublf}
For all non-empty open subset $\Omega$ of $M$:
$\exists C>0$, $\forall d>0$,
\begin{equation*}
\limsup_{T\to 0} 
T\ln\sup_{ u_{0}\in \dot{E}_{d/ T} }
\frac{\|u_{0}\|_{L^{2}(M)}}{\|e^{it\Delta}u_{0}\|_{L^{2}(\Omega_{T})}}
\leq Cd
\mbox{, where } 
\dot{E}_{d/T}=\Vect{e_{j}}_{\omega_{j}\leq d/T} \setminus \{0\}
\ .
\end{equation*}
\end{corollary}



\section{High frequencies}
\label{sec:hf}

To analyze 
how violent fast controls are  
for high frequency vibrations (cf. section~\ref{sec:plan}),
we introduce a {\em wavelength scale} 
$(h_{k})_{k\in\N}$, 
i.e. a decreasing sequence of positive real numbers converging to $0$,
and the corresponding {\em spectral scale}
$(E_{k})_{k\in\N}$ of subspaces of $L^{2}(M)$ defined by 
$E_{k}=\Vect{e_{j}}_{a<h_{k}\omj<b}$
for fixed $b>a>0$
(note that these subspaces may be overlapping),
where the spectral data $(e_{j},\omega_{j})$
are defined in (\ref{eqeigen}).

In this section, we take up a strategy of Lebeau in~\cite{Leb92}:
we reduce the observation of Schr{\"o}dinger equation on $\Omega$
in small time $hT$ 
to the observation of the semiclassical Schr{\"o}dinger equation 
$ih\partial_{t}\psi-h^{2}\Delta\psi=0$, $\psi=\psi_{0}$ at $t=0$, 
on $\Omega$ in fixed time $T>0$
(note that taking $u_{0}=\psi_{0}$, we have 
$u(t,x):=e^{-it\Delta}u_{0}=\psi(t/h,x)$).
In the first subsection, we emphasize that 
the first step of the reduction actually 
is an equivalence.
As in~\cite{Bur97}, 
we perform the semiclassical analysis with a light microlocal tool:
the  ``microlocal measures''
introduced independently by  
P.~G{\'e}rard, P.-L.~Lions and T.~Paul, and L.~Tartar, 
and first used by G.~Lebeau in control theory
(cf. \cite{Burast} for a survey).  
In the second subsection, we keep track of the controllability cost
(our presentation also shortcuts 
the estimates in the Besov space 
$\dot{B}^{0}_{2,\infty}(\mathbb{R}_{t};L^{2}(\Omega))$
in~\cite{Leb92} and~\cite{Bur97}
thanks to lemma~\ref{lem:com}).

\subsection{Semiclassical observability}

The relevant notion of observability for 
the semiclassical Schr{\"o}dinger equation is:
\begin{definition}
\label{def:scobs}
{\em Semiclassical observability} on $\Omega\Sset M $ in time $T>0$
holds when~: 
for all $\theta\in C^{\infty}\comp(\R\times M)$ such that 
$\set{\theta\neq 0}=\Omega_{T}$,
for all $b>a\geq 1/2$, 
there is an observability constant $\Csc>0$ 
and a threshold $\hbar>0$ such that:
\begin{gather*}
\forall h\in\left]0,\hbar\right], \forall \psi_{0}\in 
\Vect{e_{j}}_{a<h\omj<b},
\quad
\norm{\psi_{0}}_{L^{2}(M)}
\leq \Csc \norm{\theta e^{-ith\Delta}\psi_{0}}_{L^{2}(\R\times M)}
\ .
\end{gather*}
\end{definition}

The purpose of this subsection is to prove:
\begin{theorem}
\label{th:scobs}
Let $\Omega\Sset M$ 
and let $L_{\Omega}$ be 
the length of the longest generalized geodesic in $\adh{M}$
which does not intersect $\Omega$. 
The geodesics condition $T>L_{\Omega}$ is necessary and sufficient for
semiclassical observability on $\Omega$ in time $T$
(cf. definition~\ref{def:scobs}). 
\end{theorem}

\begin{remarks}
Semiclassical observability implies observability 
(cf. remark~\ref{rem:hfobs})
but the converse does not hold (cf. section~\ref{sec:ubprod}). 

As for the wave equation (cf.~\cite{BG97}), 
this theorem does not hold with 
the smooth characteristic function $\theta$
replaced by $\cd{1}_{\Omega_{T}}$.

If the condition $b>a\geq 1/2$ 
is replaced by $b>a\geq \aa$
then this becomes the definition of  
semiclassical observability on $\Omega$ in time $2\aa T$.
Moreover, the condition 
``~for all $\theta\in C^{\infty}\comp(\R\times M)$ such that 
$\set{\theta\neq 0}=\Omega_{T}$, for all $b>a\geq 1/2$~'' 
could be equivalently replaced by 
``~there is a $\theta\in C^{\infty}\comp(\R\times M)$ such that 
$\set{\theta\neq 0}=\Omega_{T}$, there is a $b>a= 1/2$~''.
\end{remarks}

\begin{proof}[Proof of theorem~\ref{th:scobs}]
We refer to~\cite{Bur97} and the survey~\cite{Burast} 
for the definition and properties
of semiclassical measures (a.k.a. Wigner measures)
that we use in this proof.

\medskip

We first prove the sufficiency by contradiction.
We assume $T>L_{\Omega}$ and that semiclassical observability does not hold,
i.e. there are  real numbers $b>a\geq 2^{-1/2}$,
a decreasing sequence $(h_{k})_{k\in\N}$
of positive real numbers converging to $0$, 
a sequence $(\psi_{0}^{k})_{k\in\N}$ of initial data in $L^{2}(M)$
such that: 
\begin{equation}
\label{eqnotscobs}
\forall k\in\N,\ 
\psi_{0}^{k}\in E_{k}=
\Vect{e_{j}}_{a<h_{k}\omj<b}
\ \text{and}\
\norm{\psi_{0}^{k}}_{L^{2}(M)}
> \frac{1}{k} \norm{\theta e^{-ith_{k}\Delta}\psi_{0}}_{L^{2}(\R\times M)}
\ .
\end{equation} 
We shall use the more convenient unambiguous abbreviations 
$h=h_{k}$, $\psi_{0}^{h}=\psi_{0}^{k}$ and 
$h\to 0$ instead of $k\to \infty$.
Without loss of generality, we assume $\norm{\psi_{0}^{h}}=1$
so that $\psi^{h}(t,x)=e^{-ith\Delta}\psi_{0}^{h}$ is bounded in 
$L^{2}\loc(\R\times M)$, and therefore, without loss of generality again,
we assume that $(\psi^{h})$ has a semiclassical measure $\mu$.
Note that 
$\mu(t,x,\tau,\xi)$ is a positive Radon measure on $T^{*}(\R\times M)$
which describes the asymptotic microlocal distribution 
of the space-time waves density $|\psi^{h}(t,x)|^{2}dtdx$.

The estimate in (\ref{eqnotscobs}) implies $\norm{\theta\psi^{h}}=o(1)$
so that $\mu\left(\Omega_{T}\right)=0$.
The first part of   (\ref{eqnotscobs}) says 
$\psi_{0}^{h}\in \Vect{e_{j}}_{a<h\omj<b}$
which implies $\norm{\Delta\psi_{0}^{h}}\leq (b/h)^{2}\norm{\psi_{0}^{h}}$,
hence $ih\partial_{t}\psi^{h}=h^{2}\Delta\psi^{h}$ is bounded in 
$L^{2}\loc(\R\times M)$ and, in particular, $(\psi^{h})$ is $h$-oscillating.
Therefore:
\begin{equation}
\label{eqmuI}
\text{for all non empty interval } I, \quad  
\mu\left(\set{(t,x)\in I\times M}\right)=\abs{I}>0 \ .
\end{equation}
Another consequence of $\psi_{0}^{h}\in \Vect{e_{j}}_{a<h\omj<b}$
is that $t\mapsto \psi^{h}$ is a linear combination of 
semiclassical time exponentials $t\mapsto \exp(it\tau/h)$
with $\tau\in \{-h^{2}\omega_{j}^{2}\}_{a<h\omj<b}$
so that $\supp\mu\subset\set{\tau\in \left[a^{2},b^{2}\right]}$.
From $ih\partial_{t}\psi^{h}-h^{2}\Delta\psi^{h}=0$,
it can be deduced by the symbolic calculus that 
$\supp\mu\subset\set{\tau=\abs{\xi}^{2}}$ and 
$\set{\tau-\abs{\xi}^{2},\mu}=\partial_{t}\mu-2\abs{\xi}\nabla_{x}\mu=0$.
Together with the Dirichlet boundary condition, this equation 
for $\mu$ means that, on any surface $\set{2\abs{\xi}=v}$, 
$\mu$ is invariant by the generalized geodesic flow at speed $v$.
But $\supp\mu\subset\set{\tau=\abs{\xi}^{2}\in\left[a^{2},b^{2}\right]}
\subset\set{v=2\abs{\xi}\geq 2a=1}$,
hence $\mu\left(\Omega_{T}\right)=0$ 
and the geodesics condition $T>L_{\Omega}$
imply $\mu=0$, in contradiction with (\ref{eqmuI}).

\medskip

Now we prove the necessity by contradiction.
We assume that semiclassical observability holds
and $T\leq L_{\Omega}$, i.e. there is a generalized geodesic 
$x:\left[ 0,T\right]\to \adh{M}$ which does not intersect $\Omega$.
Without loss of generality, we assume $x(0)\notin\partial M$.

To construct initial data which concentrate on  $x_{0}=x(0)$
with initial momentum $\xi_{0}$ with $|\xi_{0}|\in ]a,b[$
and with the direction corresponding to $x'(0)$,
we introduce a smooth cut-off function $\chi$ 
compactly supported in a chart of $M$ around $x_{0}$
such that $\chi=1$ in a neighborhood of $x_{0}$,
and define $\psi_{0}^{h}$ as the function 
$x\mapsto\chi(x)\exp(ix.\xi_{0}/h)\exp(-(x-x_{0})^{2}/h)$
divided by its $L^{2}(M)$ norm.
Then the semiclassical measure of $(\psi_{0}^{h})$ is 
$\delta(x-x_{0},\xi-\xi_{0})$
and, thanks to proposition~4.11 in~\cite{Bur97},
we may assume without loss of generality that 
$(\psi_{0}^{h})$ has been projected on 
$\Vect{e_{j}}_{a<h\omj<b}$.
As before, 
we may assume that $(\psi^{h})$ has a semiclassical measure $\mu$.
Taking the limit $h\to 0$ in the inequality defining 
semiclassical observability yields $\mu\left(\Omega_{T}\right)>0$. 
As before $\mu$ is invariant by the generalized geodesic flow 
at speed $v=2\abs{\xi_{0}}>2a=1$.
We may choose $\xi_{0}$ with $v$ close enough to $1$
so that the support of $\mu$ is so close to 
the image of the generalized geodesic $x$ in $T^{*}(\R\times M)$
that it does not intersect $\Omega$,
in contradiction with $\mu\left(\Omega_{T}\right)>0$. 
\end{proof}

\subsection{Upper bound at high frequencies under the geodesics condition}
\label{sec:ubhf}

The main result proved in this subsection is that 
semiclassical observability implies 
``~observability at cost $C/\sqrt{T}$ modulo low frequencies~''~:  
\begin{theorem}
\label{th:hfobs}
Semiclassical observability on $\Omega\Sset M$ in time $\TOm$
(cf. definition~\ref{def:scobs}) implies that:
$\exists \kk\in\N$, $\forall d>\TOm$, $\exists C_{d}>0$, 
$\forall k\geq \kk$, $\forall T\in[h_{k}d,h_{\kk}d]$,
\begin{gather*}
\forall u_{0}\in L^{2}(M),\
(1+O(h_{k}^{2}/T))\norm{u_{0}}_{L^{2}(M)}^{2}\leq 
\frac{C_{d}^{2}}{T}\norm{e^{-it\Delta}u_{0}}_{L^{2}(\Omega_{T})}^{2}
+\norm{\pi_{k}u_{0}}_{L^{2}(M)}^{2}
\ ,
\end{gather*}
where $h_{k}=2^{-k}$ is the dyadic scale
and $\pi_{k}$ is the projection on 
$\Vect{e_{j}}_{h_{k}\omj\leq1}$.
\end{theorem}

\begin{remark}
\label{rem:hfobs}
As usual, since $\pi_{k}$ is a compact operator, 
we can get rid of the remainder low frequency term 
in the observability inequality of theorem~\ref{th:hfobs}
by the unique continuation property of elliptic operators
(as in lemma~6 in~\cite{Leb92}).
Hence, under the geodesics condition, 
theorems~\ref{th:scobs} and~\ref{th:hfobs} imply 
the exact controllability of Schr{\"o}dinger equation 
from $\Omega$ in any time, 
i.e. $C_{T,\Omega}<+\infty$ for all $T>0$
(this is the analogue for interior controllability 
of the boundary controllability theorem in~\cite{Leb92}).
\end{remark}

Taking $T=h_{k}d$ in this theorem and combining it with theorem~\ref{th:scobs}
allow us to state a counterpart to 
the upper bound at low frequencies of corollary~\ref{cor:ublf},
i.e. the following upper bound at high frequencies: 
\begin{corollary}
\label{cor:ubhf}
Let $\Omega\Sset M$ 
and let $L_{\Omega}$ be 
the length of the longest generalized geodesic in $\adh{M}$
which does not intersect $\Omega$. 
For all $d>L_{\Omega}$, 
there is a constant $C_{d}>0$ 
and a dyadic sequence of positive times $T$ converging to $0$
such that:
\begin{gather*}
\forall u_{0}\in 
\Vect{e_{j}}_{\omj\geq d/T},
\quad
\norm{u_{0}}_{L^{2}(M)}
\leq \frac{C_{d} }{\sqrt{T}}
\norm{ e^{-it\Delta}u_{0} }_{L^{2}(\Omega_{T})}
\ .
\end{gather*}
\end{corollary}
The dual statement is a smoothing property at low control cost:
\begin{corollary}
\label{cor:smoothing}
For all $d>L_{\Omega}$, 
there is a constant $C_{d}>0$ 
and a dyadic sequence of positive times $T$ converging to $0$
such that:
$\forall u_{0}\in L^{2}(M)$, 
$\exists g\in L^{2}(\mathbb{R}\times M)$
such that $\|g\|_{L^{2}(\mathbb{R}\times M)}
\leq \|\Pi_{d/T} u_{0}\|_{L^{2}(M)}C_{d}/\sqrt{T}$
and the solution $u\in C^{0}([0,\infty);L^{2}(M))$
of (\ref{eqSchr}) with Cauchy data
$u=u_{0}$ at $t=0$,
satisfies 
$u=u_{T}\in E_{d/T}\subset C^{\infty}({\adh{M}})$ at $t=T$,
where $\Pi_{d/T}$
is the projection on the space orthogonal to 
$E_{d/T}=\Vect{e_{j}}_{\omega_{j}< d/T}$.
In particular: $\forall s\geq 0$,
$\|u_{T}\|_{H^{s}(M)}
\leq (1+C_{d})(d/T)^{s}\|u_{0}\|_{L^{2}(M)}$.
\end{corollary}

The first preliminary step in proving theorem~\ref{th:hfobs}
is to deduce a ``~high frequency observability inequality~'' 
from semiclassical observability:
\begin{lemma}
\label{lem:hfobs}
Semiclassical observability on $\Omega\Sset M$ in time $T>0$
with $b=2=a^{-1}$ implies that
there is an observability constant $\Chf>0$ 
and a threshold $\kk\in\N^{*}$ such that:
$\forall k\geq \kk$, $\forall S\in[h_{k}T,h_{\kk-1}T]$,
\begin{gather*}
\forall v_{0}\in E_{k}=\Vect{e_{j}}_{h_{k-1}^{-1}<\omj<h_{k+1}^{-1}},
\quad
\norm{v_{0}}_{L^{2}(M)}
\leq \frac{\Chf}{\sqrt{S}}\norm{e^{-it\Delta}v_{0}}_{L^{2}(\Omega_{S})}
\ .
\end{gather*}
\end{lemma}
\begin{proof}
We choose $\tOm\in C^{\infty}\comp(M)$ and $\tT\in C^{\infty}\comp(\mathbb{R})$
with values in $[0,1]$ such that 
$\set{\tOm\neq 0}=\Omega$,
and $\set{\tT\neq 0}=\left]0,T\right[$.
Let $\Csc$ and $\hbar$ be the positive constants 
obtained by applying definition~\ref{def:scobs} 
with $\theta(t,x)=\tT(t)\tOm(x)$ and $b=2=a^{-1}$.
Choosing $\kk$ such that $h_{\kk}<\hbar$, 
the semiclassical observability inequality of definition~\ref{def:scobs} 
implies:
\begin{gather*}
\forall k\geq \kk, \forall \psi_{0}\in E_{k}, 
\quad
\norm{\psi_{0}}_{L^{2}(M)}^{2}
\leq \Csc^{2}\int\norm{\tOm e^{-ish_{k}\Delta}\psi_{0}}_{L^{2}(M)}^{2}
\abs{\tT(s)}^{2}ds
\ .
\end{gather*}
The change of variable $t=sh_{k}$ and the definition of $\tOm$ and $\tT$ yield:
\begin{gather*}
\norm{\psi_{0}}_{L^{2}(M)}^{2} \leq 
\frac{\Csc^{2}}{h_{k}}
\int_{0}^{h_{k}T}\norm{e^{-it\Delta}\psi_{0}}_{L^{2}(\Omega)}^{2}dt
\ .
\end{gather*}
Taking $\psi_{0}=e^{-iNh_{k}T}v_{0}$ and changing $t$ by a translation yields:
\begin{gather*}
\forall N\in\N, \forall k\geq \kk, \forall v_{0}\in E_{k}, 
\
\norm{v_{0}}_{L^{2}(M)}^{2} \leq 
\frac{\Csc^{2}}{h_{k}}
\int_{Nh_{k}T}^{(N+1)h_{k}T}\norm{e^{-it\Delta}v_{0}}_{L^{2}(\Omega)}^{2}dt
\ .
\end{gather*}
Let $k\geq m\geq \kk$. Summing up from $N=0$ to $N=h_{m}h_{k}^{-1}-1$,
multiplying by $h_{k}h_{m}^{-1}$ and setting $\Chf=\Csc\sqrt{2T}$ yield:
\begin{gather*}
\forall v_{0}\in E_{k}, 
\
\norm{v_{0}}_{L^{2}(M)}^{2} \leq 
\frac{\Chf^{2}}{h_{m-1}T}
\int_{0}^{h_{m}T}\norm{e^{-it\Delta}v_{0}}_{L^{2}(\Omega)}^{2}dt
\ .
\end{gather*}
This inequality completes the proof of lemma~\ref{lem:hfobs}
since, for all $k\geq \kk$ and $S\in[h_{k}T,h_{\kk-1}T]$,
there is a $m\in[\kk,k]$ such that $S\in[h_{m}T,h_{m-1}T]$.
\end{proof}

The second preliminary step in proving theorem~\ref{th:hfobs}
is to introduce a time frequency decomposition which is 
semiclassically equivalent to 
the spatial decomposition into the spectral scale $(E_{k})$.
This provides an easy way to overcome
the following difficulty
(cf. lemma~\ref{lem:com}): 
multiplication by $\theta$ 
(which corresponds to observing on $]0,T[\times\Omega$) 
does not commute with the projection on $E_{k}$.
(Note that this difficulty is even greater in boundary observability
but can be overcome by the analogue of lemma~\ref{lem:com}).

The Fourier transform of $v\in L^{2}(\R\times M)$ with respect to $t$ is~:
\begin{gather*}
\hat{v}(\tau,x)=\int e^{-i\tau t }v(t,x)dt \, .
\end{gather*}
For any $\phi\in L^{\infty}(\mathbb{R})$,
the frequency cut-off $\phi(D_{t})$
and the spectral cut-off $\phi(\sD)$
are the bounded operators on $L^{2}$ defined by:
\begin{align*}
&\forall v\in L^{2}(\R\times M), 
&\left(\phi(D_{t})v\right)(t,x)
&=\frac{1}{2\pi}\int e^{i\tau t}\phi(\tau)\hat{v}(\tau,x)d\tau 
\ ,
\\ 
&\forall v\in L^{2}(M),
&\phi(\sD)v
&=\sum_{j\in \N^{*}}\phi(\omega_{j})(v|e_{j})_{L^{2}(M)}e_{j}
\ .
\end{align*}
For instance, 
the projection on $E_{k}$ writes 
$\cd{1}_{[h_{k-1}^{-1},h_{k+1}^{-1}]}(\sD)=\cd{1}_{[1/2,2]}(h_{k}\sD)$
with this notation.
For $\phi \in \mathcal{S}(\mathbb{R})$,
these cut-off operators extend to $L^{2}(M, \mathcal{S}^{'}(\mathbb{R}_{t}))$
and satisfy the following ``~compatibility~'' relations~:
\begin{gather*}
\phi(D_{t})\left(e^{it\Delta}e_{j}\right)
=\phi(\omega_{j})e^{it\omega_{j}^{2}}e_{j}
=e^{it\Delta}\phi(\sD)e_{j}
\ .
\end{gather*}
We shall need the following commutator estimate:
\begin{lemma}
\label{lem:com}
For all $\theta_{T}\in C^{\infty}\comp(\left] 0,T\right[)$,
$\phi \in \mathcal{S}(\mathbb{R})$, and $v\in L^{2}(\R\times M)$~:
\begin{gather*}
\norm{[\, \theta_{T},\phi(D_{t}) \, ]Wv}_{L^{2}(\R\times M)}
\leq 
\norm{v}_{L^{2}(\R\times M)}
(1+T)\, \norm{\partial_{t}\theta_{T}}_{L^{\infty}}
\int (1+\abs{t})\abs{t\hat{\phi}(t)}\, dt
 \, ,
\end{gather*}
where $W$ denotes the weight multiplication $Wv(t,x)=(1+|t|) v(t,x)$ 
and $[\, \theta_{T},\phi(D_{t}) \, ]$
denotes the commutator $\theta_{T}\phi(D_{t})-\phi(D_{t})\theta_{T}$ 
between the multiplication by $\tT$
and the time frequency cut-off $\phi(D_{t})$.
\end{lemma}
\begin{proof}
Since the operator does not act in the $x$ variable, 
we may forget about $x$ and write its kernel as~:
\begin{gather*}
K(t,s)=(1+|s|)\left( \theta_{T}(t)-\theta_{T}(s) \right) \hat{\phi}(s-t)
\, . \end{gather*}
By Schur's lemma, the bound sought for this operator 
will result from the same bound on 
$\sup_{t} \int|K(t,s)|\, ds$
and on 
$\sup_{s} \int|K(t,s)|\, dt$.

Using $\supp \theta_{T} \subset\, ]0,T[$
and Taylor's inequality yields~:
\begin{multline*}
\int \abs{K(t,s)}\, ds =  
\int (1+\abs{s+t})\, \abs{\theta_{T}(t)-\theta_{T}(s+t)}\, 
\abs{\hat{\phi}(s)} \, ds \\
\leq  
\int (1+T)(1+|t|) \, \abs{s} \, \norm{\partial_{t}\theta_{T}}_{L^{\infty}}\, 
\abs{\hat{\phi}(s)} \, ds
\, ,
\end{multline*}
and similarly $\int|K(t,s)|\, dt\leq  
\int (1+T)(1+\abs{t})\, \abs{t}\, \norm{\partial_{t}\theta_{T}}_{L^{\infty}}
\abs{\hat{\phi}(t)} \, dt$.
\end{proof}

We choose a smooth real valued cut-off function $\tau\mapsto \chi(\tau)$
the time frequency parameter $\tau\in\R$, such that:
\begin{gather}
\label{eqchisupp}
\supp\chi\subset  \set{1/2<\abs{\tau}<2} \ ,
\\
\label{eqchisum}
\sum_{k\geq 0}\chi^{2}(2^{-k}\tau)\geq 1 \text{ on } \set{\abs{\tau}>1} \ .
\end{gather}
For any $m\in\N^{*}$, 
since $h_{m}|\omega_{j}|>1$ and $h_{k}|\omega_{j}|<2$ imply $k\geq m$,
(\ref{eqchisupp}) and (\ref{eqchisum}) imply:
\begin{equation}
\label{eqchi1}
\forall v\in L^{2}(M), \norm{v-\pi_{m}v}_{L^{2}(M)}^{2}
\leq \sum_{k\geq m} \norm{\chi(h_{k}\sD)v}_{L^{2}(M)}^{2}
\ .
\end{equation}  
Setting $C_{\chi}=2\norm{\chi}_{L{\infty}}^{2}$, (\ref{eqchisupp}) implies
$\sum_{k\in\N}\chi^{2}(2^{-k}\tau)\leq C_{\chi}$
(for each $\tau$ there are at most two nonzero terms in this sum),
so that:
\begin{equation}
\label{eqchi2}
\forall v\in L^{2}(M\times M), 
\sum_{k\in\N}\norm{\chi(h_{k}D_{t})v}_{L^{2}(\R\times M)}^{2}
\leq C_{\chi}\norm{v}_{L^{2}(\R\times M)}^{2}
\ .
\end{equation}  

\begin{proof}[Proof of theorem~\ref{th:hfobs}]
We assume semiclassical observability on $\Omega\Sset M$ in time $\TOm$
and apply lemma~\ref{lem:hfobs}.
Let $d>\TOm$, $k\geq \kk$, and $T\in[h_{k}d,h_{\kk}d]$.
Applying the high frequency observability inequality of lemma~\ref{lem:hfobs}
with $S=T/2$ and $v_{0}=e^{-iS\Delta}\chi(h_{m}\sD)u_{0}$,
choosing $\tT\in C^{\infty}\comp(\mathbb{R})$
with values in $[0,1]$ such that $\set{\tT\neq 0}=\left]0,T\right[$
and $\set{\tT=1}=\left]S/2,3S/2\right[$, 
and setting $C_{d}'=\sqrt{d/\TOm}\Chf$, we obtain:
\begin{multline}
\label{eqhf}
\forall m\geq k,\forall u_{0}\in L^{2}(M),
 \norm{\chi(h_{m}\sD)u_{0}}_{L^{2}(M)} \\
\leq  \frac{\Chf}{\sqrt{S}} 
\norm{\chi(h_{m}\sD)u}_{L^{2}(\left]\frac{S}{2},\frac{3S}{2}\right[\times\Omega)}
\leq  \frac{C_{d}'}{\sqrt{T}} 
\norm{\tT\chi(h_{m}\sD)u}_{L^{2}(\R\times\Omega)}
\ .
\end{multline}
Since $t\mapsto e^{-it\Delta}$ is unitary: 
\begin{gather*}
\int_{\R\times M}\abs{(1+\abs{t}^{-1})u(t,x)}^{2}dxdt\leq 
\int_{M}\abs{u(t,x)}^{2}dx 
\int_{\R}\abs{1+\abs{t}^{-1}}^{2}dt  
\ .
\end{gather*}
Applying lemma~\ref{lem:com} with $\psi{\tau}=\chi(h_{m}\tau)$ 
and $v(t,x)=(1+\abs{t}^{-1})\cd{1}_{\Omega}(x)u(t,x)$,
and setting $C=\norm{1+\abs{t}^{-1}}_{L^{2}}
(1+T)\, \norm{\partial_{t}\theta_{T}}_{L^{\infty}}
\int (1+\abs{t})\abs{t\hat{\phi}(t)}\, dt$,
we obtain:
\begin{gather*}
\norm{\tT\chi(h_{m}\sD)u-\chi(h_{m}D_{t})\tT u}_{L^{2}(\R\times \Omega)}^{2}
\leq C^{2}h_{m}^{2}\norm{u_{0}}_{L^{2}(\Omega)}^{2}
\ .
\end{gather*}
Since $\abs{\tT}\leq 1$, this combines with (\ref{eqhf}) into:
$\forall m\geq k$, $\forall u_{0}\in L^{2}(M)$, 
\begin{gather*}
\norm{\chi(h_{m}\sD)u_{0}}_{L^{2}(M)}^{2}
\leq \frac{2C_{d}'{}^{2}}{T}
\left(
\norm{\chi(h_{m}D_{t})u}_{L^{2}(\Omega_{T})}^{2}
+C^{2}h_{m}^{2}\norm{u_{0}}_{L^{2}(\Omega)}^{2}
\right)
\ . 
\end{gather*}
Summing up this inequality for $m\geq k$,
we obtain,  
thanks to (\ref{eqchi1}), (\ref{eqchi2}) and 
$\sum_{m\geq k}h_{m}^{2}=4 h_{k}^{2}/3$~:
\begin{gather*}
\norm{u_{0}-\pi_{k}u_{0}}_{L^{2}(M)}^{2}\leq 
\frac{2C_{d}{}^{2}}{T}
\left(C_{\chi}
\norm{u}_{L^{2}(\Omega_{T})}^{2}
+\frac{4}{3}C^{2}h_{k}^{2}\norm{u_{0}}_{L^{2}(\Omega)}^{2}
\right)
\ .
\end{gather*} 
Adding $\norm{\pi_{k}u_{0}}_{L^{2}(M)}^{2}
-\frac{8}{3}C_{d}'{}^{2}C^{2}h_{k}^{2}\norm{u_{0}}_{L^{2}(\Omega)}^{2}$
on both sides 
completes the proof of theorem~\ref{th:hfobs}
with  $C_{d}=\sqrt{2d/\TOm}\Chf C_{\chi}$. 
\end{proof}


\section{A window problem for nonharmonic Fourier series}
\label{sec:nonharm}

In this section we prove theorem~\ref{th:1dSL}
which generalizes theorem~\ref{th:1d} to Sturm-Liouville operators
(in particular
to a segment with any Riemannian metric).
By spectral analysis, 
it reduces to a refinement (cf. proposition~\ref{prop:windpb})
of the well studied 
window problem for nonharmonic Fourier series (cf.~\cite{SAI00})
which is solved in subsections~\ref{sec:windpb} and~\ref{sec:entfunc}
following \cite{LMheatcost} quite closely. 

\subsection{Boundary control of a segment}
\label{sec:1d}

Let $X>0$. 
We consider the Sturm-Liouville operator $A$ on $L^{2}(0,X)$ 
with domain $D(A)$ defined by 
\begin{align*}
(Af)(x) &= \left(p(x)f'(x)\right)' +q(x)f(x)
\quad \text{ for } x\in [0,X]
\\
D(A) &= \{f\in H^{2}(0,X) \,|\,
\left(a_{0}f+b_{0}f'\right)(0) = 0
=\left(a_{1}f+b_{1}f'\right)(X)
\}
\end{align*}
where all the coefficients are real and satisfy:
\begin{gather}
\label{eqSLAssum}
a_{0}^{2}+b_{0}^{2}= a_{1}^{2}+b_{1}^{2}=1 \, ,\quad
0<p\in C^{2}([0,X]) 
\, ,\quad
q\in C^{0}([0,X])
\ ,
\end{gather}
Under these assumptions, 
$-A$ is self-adjoint
and has a sequence $\{\lambda_{n}\}_{n\in \mathbb{N}^{*}}$
of increasing eigenvalues 
and an orthonormal Hilbert basis $\{e_{n}\}_{n\in \mathbb{N}^{*}}$ 
in $L^{2}(0,X)$ of corresponding eigenfunctions, i.e.: 
\begin{equation*}
\forall n\in \mathbb{N}^{*}, \quad -Ae_{n}  = \lambda_{n}e_{n} 
\quad \mbox{and} \quad 
\lambda_{n}<\lambda_{n+1} \ .
\end{equation*}
Moreover, $(\ref{eqSLAssum})$ ensures the following eigenvalues asymptotics:
\begin{equation}
\label{eqSLEig}
\exists \nu\in \mathbb{R},\
\lambda_{n}=\frac{\pi^{2}}{L^{2}}\left(n+\nu\right)^{2} + O(1) 
\ \mbox{as} \  n\to \infty \ ,
\quad \mbox{where} \quad L=\int_{0}^{X}\sqrt{p(x)}\, dx
\ .
\end{equation}
We use the following notations for the Sobolev spaces based on $A$:
\begin{gather*}
H^{0}_{A}(0,X)=L^{2}(0,X)
\quad\text{ and }\quad
H^{1}_{A}(0,X)=\adh{D(A)}^{H^{1}}
\ .
\end{gather*}

\begin{theorem}
\label{th:1dSL}
For any $\alpha > \alpha^{*}$ defined by $(\ref{eqalphastar})$,
there exists $C>0$ such that, 
for any coefficients $(\ref{eqSLAssum})$, 
setting $k=1$ if $b_{1}=0$ and $k=0$ otherwise,
for all $T\in\, ]0,\inf(\pi,L)^{2}]$
and $u_{0}\in H^{k}_{A}(0,X)$
the solution $u\in C^{0}([0,\infty);H^{k}_{A}(0,X))$
of 
\begin{gather*} 
i\partial_{t}u= \partial_{x}\left(p(x)\partial_{x}u\right) +q(x)u
\quad \text{for}\ (t,x)\in\, ]0,T[\times ]0,X[\ , 
\\ 
\left(a_{0}u+b_{0}\partial_{x}u\right)\res{x=0} = 0 
=\left(a_{1}u+b_{1}\partial_{x}u\right)\res{x=X} 
\quad\text{and}\quad
u\res{t=0} = u_{0} \ ,
\end{gather*}
satisfies 
$\displaystyle
\|u_{0}\|_{H^{k}_{A}(0,X)} 
\leq  C\exp(\alpha L^{2}/T)
\|\partial_{x}^{k}u\res{x=X}\|_{L^{2}(0,T)}
$.
\end{theorem}

\subsection{Reduction of the window problem to a problem on entire functions}
\label{sec:windpb}

First note that the theorem~\ref{th:1dSL} can be reduced to the case 
$\lambda_{1}>0$ 
by the multiplier $t\mapsto \exp(i\lambda t)$,
to the case $L=\pi$ 
by the time rescaling $t\mapsto \sigma t$
with $\sigma=(\pi/L)^{2}$,
and to the time interval $[-T/2,T/2]$
by the time translation $t\mapsto t-T/2$.

From now on  we assume
$\lambda_{1}>0$ and $L=\pi$.
Making a weaker assumption on the remainder term in $(\ref{eqSLEig})$,
we shall only use the following spectral assumption:
\begin{equation}
\label{eqSpec}
\forall n\in \mathbb{N}^{*},\  
0<\lambda_{n}<\lambda_{n+1} 
\quad \mbox{and} \quad 
\lambda_{n}=n^{2} + O(n)
\ \mbox{as} \ n\to \infty \ .
\end{equation}

In terms of the coordinates
$c=(c_{k})_{k\in \mathbb{N}^{*}}$ of $A^{k/2}u_{0}$ 
in the Hilbert basis $(e_{k})_{k\in \mathbb{N}^{*}}$,
we have to solve the following window problem:

\begin{proposition}
\label{prop:windpb}
For any $\alpha > \alpha^{*}$,
there exists $C>0$ such that, 
for all $(\lambda_{n})_{n\in \mathbb{N}^{*}}$ satisfying (\ref{eqSpec}),
for all $T\in \left]0,\pi\right]$: 
\begin{equation}
\label{eq:windpb}
\forall c\in l^{2}(\N^{*}), \quad
\norm{c}_{l^{2}}
\leq  Ce^{\alpha\pi^{2}/T }
\|f\|_{L^{2}(-\frac{T}{2},\frac{T}{2})}
\text{ where } f(t)= \sum_{n=1}^{\infty} c_{n}e^{i\ldn t}
\ .
\end{equation}
\end{proposition}
                                
The well-known method to study the nonharmonic Fourier series $f$
is to construct a sequence $(g_{n})_{n\in \mathbb{N}^{*}}$ in $L^{2}(-T/2,T/2)$
which is bi-orthogonal to the sequence 
$\{\exp(-\lambda_{n}t)\}_{n\in \mathbb{N}^{*}}$, i.e.
\begin{equation}
\label{eqbiorth}
\int_{-T/2}^{T/2}g_{n}(t)e^{-\ldn t} \, dt=1 
\quad \mbox{ and  } \quad
\forall k\in \mathbb{N}^{*} ,\, k\neq n ,\,  
\int_{-T/2}^{T/2}g_{n}(t)e^{-\ldk t} \, dt=0  \ . 
\end{equation}
Then: $\norm{c}_{l^{2}}^{2}=\sum_{n}(f,g_{n})\conj{c_{n}}
=(f,\sum_{n}g_{n}c_{n})
\leq \norm{f}_{L^{2}}\norm{\sum_{n}g_{n}c_{n}}_{L^{2}}$.
Hence, introducing the Gramm operator $G$ on $l^{2}(\N^{*})$
defined by the coefficients $(g_{n},g_{k})_{L^{2}}$
for $n$ and $k$ in $\N^{*}$,
(\ref{eq:windpb}) results from $\norm{G}\leq Ce^{\alpha/T }$.
But, applying Schur's lemma to the self-adjoint operator $G$ yields:
$\norm{G}^{2}\leq \sup_{n}\sum_{k}\abs{(g_{n},g_{k})_{L^{2}}}$.
Thus, to prove proposition~\ref{prop:windpb}
it is enough to construct bi-orthogonal functions $g_{n}$
with good growth estimates of their scalar products as $T$ tends to zero.
We shall readily explain how  
the following proposition on entire functions
yields this construction,
and postpone its proof to subsection~\ref{sec:entfunc}.
\begin{proposition}
  \label{prop:G}
Let $\alpha^{*}$ be defined by $(\ref{eqalphastar})$.
Let $\{\lambda_{n}\}_{n\in \mathbb{N}^{*}}$ be a sequence of real numbers 
satisfying $(\ref{eqSpec})$.
For all $\eps>0$
there is a $C_{\eps}>0$ such that, 
for all $\tau\in ]0,1]$,
there is a sequence of entire functions $\{G_{n}\}_{n\in \mathbb{N}^{*}}$ 
satisfying, 
for all $n$ and $k$ in $\mathbb{N}^{*}$: 
\begin{align}
&
  G_{n} \mbox{ is of exponential type }\tau
\mbox{, i.e. } \limsup_{r\to +\infty} r^{-1}\sup_{|z|=r}\ln |G_{n}(z)| \leq\tau,
\label{eqG1} \\
&
G_{n}(\ldn)=1 
\quad \text{ and  } \quad
G_{n}(\ldk)=0 \ \text{ if  } \ k\neq n, 
\label{eqG2} \\
&
\abs{ (G_{n},G_{k})_{L^{2}} }\leq C_{\eps} 
e^{ -\eps\sqrt{ \abs{\ldn-\ldk}/2} } e^{\alpha^{*}(\pi+\sqrt{2}\eps)^{2}/\tau}
\ .
\label{eqG3} 
\end{align}
\end{proposition}

According to the Paley-Wiener theorem (1934),
$(\ref{eqG1})$ implies that the function $x\mapsto G_{n}(x)$
is the unitary Fourier transform of a function $t\mapsto g_{n}(t)$
in $L^{2}(\mathbb{R})$ supported in $[-\tau,\tau]$.
With $\tau=T/2$, this yields:
\begin{equation}
\label{eqg}
G_{n}(x)=\frac{1}{\sqrt{2\pi}}\int_{-T/2}^{T/2}g_{n}(t)e^{-itx}\, dt
\quad \mbox{and} \quad 
\|g_{n}\|_{L^{2}}=\|G_{n}\|_{L^{2}} \ .
\end{equation}
Hence $(\ref{eqG2})$ implies $(\ref{eqbiorth})$   
and $(\ref{eqG3})$ implies that:
\begin{gather*}
\norm{G}^{2}\leq \sup_{n}\sum_{k}\abs{(g_{n},g_{k})_{L^{2}}}
\leq C_{\eps} e^{2\alpha^{*}(\pi+\sqrt{2}\eps)^{2}/T}
\sup_{n}\left(1+\sum_{k\neq n }
e^{ -\frac{\eps}{\sqrt{2}}\sqrt{ \abs{\ldn-\ldk}} }
\right)
\ .
\end{gather*}

To complete the proof of proposition~\ref{prop:windpb},
we just have to estimate the last sum uniformly with respect to $n$.
For this purpose, we introduce the counting function of the sequence 
$(|\ldk-\ldn|)_{k\in \mathbb{N}^{*}\setminus\{n\}}$
for every $n\in \mathbb{N}^{*}$:
\begin{gather*}
N_{n}(r)=\#\{ k\in \mathbb{N}^{*}  \ | \ 0<\abs{\ldk-\ldn}\leq r \}
\ .
\end{gather*}
From the spectral asymptotics (\ref{eqSpec}),
we deduce that:
\begin{equation}
\label{eqN}
\exists \rr >0, \forall r\in ]0, \rr[, N_{n}(r)=0 \ , \quad  
\exists A>0, \ \forall r, \   |\sqrt{r}-N_{n}(r)|\leq A 
\ .
\end{equation}
Indeed, (\ref{eqSpec}) implies $\abs{\sldn-n}\leq C$ for some $C>0$.
If $\ldn\leq r$ then $\sqrt{\ldn+r}\leq \sqrt{2r}$ so that 
$N_{n}(r)\leq \sqrt{2r}+C$.
If $\ldn>r$ then $N_{n}(r)\leq \sqrt{\ldn+r}-\sqrt{\ldn-r}+2C\leq \sqrt{2r}+2C$.
Now the last sum writes:
\begin{multline*}
\sum_{k\neq n}
e^{ -\frac{\eps}{\sqrt{2}}\sqrt{ \abs{\ldn-\ldk}} }
=\int_{0}^{\infty}
\exp\left( -\frac{\eps}{\sqrt{2}}\sqrt{r} \right)\, dN_{n}(r) \\
=\frac{\eps}{2\sqrt{2}}\int_{0}^{\infty}\frac{N_{n}(r)}{\sqrt{r}}
\exp\left( -\frac{\eps}{\sqrt{2}}\sqrt{r} \right)\, dr
\leq \eps\int_{0}^{\infty}(2s+A)e^{-\eps s}\, ds \ .
\end{multline*}
This completes the proof that proposition~\ref{prop:G}
implies proposition~\ref{prop:windpb}
which implies theorem~\ref{th:1dSL}.

\subsection{Entire functions construction}
\label{sec:entfunc}

In this subsection, we prove proposition~\ref{prop:G}.
We follow a classical method in complex analysis:
for all $n\in \mathbb{N}^{*}$ and small $\tau>0$, 
we shall form, in a first lemma, 
an infinite product $F_{n}$
normalized by $F_{n}(\ldn)=1$
with zeros at $\ldk$ for every positive integer $k\neq n$,
and construct, in a second lemma, 
a multiplier $M_{n}$ of exponential type $\tau$
with fast decay at infinity on the real axis 
so that $G_{n}=M_{n}F_{n}$ is in $L^{2}$ on the real axis.

\begin{lemma}
\label{lemF}
Let $\{\lambda_{n}\}_{n\in \mathbb{N}^{*}}$ be a sequence of real numbers 
satisfying $(\ref{eqSpec})$.
For all $\eps>0$
there is a $A_{\eps}>0$ such that, 
for all $n\in \mathbb{N}^{*}$, 
the entire function $F_{n}$
defined by 
$\displaystyle
F_{n}(z)=\prod_{k\neq n} \left( 1-\frac{z-\ldn}{\ldk-\ldn}\right)
$ 
satisfies
\begin{equation}
  \ln|F_{n}(z-\ldn)| \leq (\sqrt{2}\pi+\eps)\sqrt{|z|} + A_{\eps}
\label{eqf1} 
\end{equation}
\end{lemma}
\begin{proof}
To prove $(\ref{eqf1})$,
we estimate the left hand side in terms of $N_{n}$:
\begin{multline*}
\ln|F_{n}(z+\ldn)| \leq  
\sum_{k\neq n} \ln\left( 1+ \frac{|z|}{|\ldk-\ldn|} \right)
=\int_{0}^{\infty} \ln\left( 1+ \frac{|z|}{r} \right) dN_{n}(r)
\\
=\int_{0}^{\infty} N_{n}(r)\frac{|z|}{|z|+r} \frac{dr}{r}
=\int_{0}^{\infty} \frac{N_{n}(|z|s)}{1+s} \frac{ds}{s}
\end{multline*}
To estimate this last integral 
we use $(\ref{eqN})$ and the integral computations:
\begin{align*}
\int_{0}^{\infty} \frac{\sqrt{s}}{1+s} \frac{ds}{s}
=\int_{0}^{\infty} \frac{2 dr}{1+r^{2}}=\pi 
\ , \quad
\int_{\frac{\rr}{|z|}}^{\infty} 
\frac{ds}{s(1+s)}
=\left[ \ln\left|\frac{s}{1+s}\right|\right]_{\frac{\rr}{|z|}}^{\infty} 
=\ln(1+\frac{|z|}{\rr})
\end{align*}
Thus we obtain 
$\ln|f_{n}(z)|\leq 
\pi\sqrt{2|z|} + A\ln(1+\frac{|z|}{\rr})$,
so that, for all $\eps>0$ there is a $A_{\eps}>0$ such that   
$ \ln|f_{n}(z)|\leq (\sqrt{2}\pi+\eps)\sqrt{|z|} + A_{\eps}$.
\end{proof}

We quote the following lemma from~\cite{LMheatcost}:
\begin{lemma}
\label{lemM}
Let $\alpha^{*}$ be defined by $(\ref{eqalphastar})$.
For all $d>0$ there is a $D>0$
such that for all $\tau>0$, 
there is an even entire function $M$ of exponential type 
(lower or equal to) $\tau$ satisfying: 
$M(0)=1$ and
\begin{equation}
\label{eqM}
\forall x>0, \quad \ln|M(x)|\leq \frac{\alpha^{*}d^{2}}{4\tau}+D-d\sqrt{x} 
\ .
\end{equation}
\end{lemma}

To prove proposition~\ref{prop:G},
we use lemmas~\ref{lemF} and~\ref{lemM}
with $d=\sqrt{2}\pi+2\eps$ and define:
$G_{n}=F_{n}M_{n}$ with $M_{n}(z)=M(z-\ldn)$.
The entire function $G_{n}$ 
has the same exponential type as $M$
since $(\ref{eqf1})$ 
implies that the exponential type of $F_{n}$ is $0$.
Hence $(\ref{eqG1})$ holds.
Equation $(\ref{eqG2})$ is an obvious consequence of 
$M_{n}(\ldn)=M(0)=1$ and the definition of $F_{n}$. 
Since $d=\pi+2\eps$ and $M$ is even, 
$(\ref{eqf1})$ and $(\ref{eqM})$ imply
\begin{gather*}
 \forall x\in \mathbb{R}, \quad 
\ln|G_{n}(x+\ldn)| \leq 
D_{\eps}+A_{\eps}-\eps\sqrt{|x|} + \frac{\alpha^{*}d^{2}}{4\tau} 
\ .
\end{gather*}
Setting $\displaystyle
C_{\eps}=2e^{2(D_{\eps}+A_{\eps})}
\int_{0}^{+\infty} e^{ -\eps\sqrt{s} } \, ds$
and $\Delta= \abs{\ldn-\ldk}/2$, 
this yields $(\ref{eqG3})$:
\begin{multline*}
\abs{ (G_{n},G_{k})_{L^{2}} }
\leq e^{2(D_{\eps}+A_{\eps}) + \frac{\alpha^{*}d^{2}}{2\tau} }
\int_{-\infty}^{+\infty} 
e^{-\eps\sqrt{\abs{x+\ldn}}-\eps\sqrt{\abs{x-\ldn}} }
\, dx\\
\leq e^{2(D_{\eps}+A_{\eps}) + \frac{\alpha^{*}d^{2}}{2\tau} }
\int_{-\infty}^{+\infty} 
e^{-\eps\sqrt{\abs{s+\Delta}}-\eps\sqrt{\abs{s-\Delta}} }
\, ds
\leq
C_{\eps} e^{ -\eps\sqrt{ \Delta } }e^{\alpha^{*}(\pi+\sqrt{2}\eps)^{2}/\tau}
\end{multline*}
Thus proposition~\ref{prop:G} is proved.

\begin{remarks}
Under the assumption $(\ref{eqSpec})$,
lemma~3 in \cite{SAI00} (which applies to more general sequences)
proves that the function 
$\displaystyle
F_{n}(z)=\prod_{k\neq n} \left[ 1- 
\left( \frac{z-\ldn}{\ldk-\ldn} \right)^{2}
\right]
$ 
satisfies $\ln|F_{n}(\ldn + z)| \leq 2\pi\sqrt{|z|}$.
In lemma~\ref{lemF}, 
the constant $2\pi$ improves to $\sqrt{2}\pi$.
We do not know if the optimal constant is $\pi$ as in
lemma~4.3 in~\cite{LMheatcost}.

Seidman obtained lemma~\ref{lemM} for
$\alpha^{*}=2\beta^{*}$
with $\beta^{*}\approx 42.86$
in the proof of Theorem~3.1 in \cite{Sei84}.
His later Theorem~1 in \cite{Sei86}
improves the rate to $\alpha^{*}=4\beta^{*}$ 
with $\beta^{*}\approx 4.17$.
Theorem~2 in \cite{SAI00}, 
which applies to much more general spectral sequences,
yields lemma~\ref{lemM} for $\alpha^{*}=48$.
As explained in ~\cite{LMheatcost}, 
lemma~\ref{lemM} does not hold for $\alpha^{*}<1/2$
and it is an interesting problem of entire function analysis 
to determine the smallest value of $\alpha^{*}$ for which it does.
\end{remarks}


\section{Upper bound under the geodesics condition}
\label{sec:ub}

In this section we prove 
theorem~\ref{th:ub}.
$ \mathcal{D}'(\mathcal{O})$ denotes the space of 
distributions on the open set $\mathcal{O}$
endowed with the weak topology 
and $\mathcal{M}(\mathcal{O})$ denotes the subspace of 
Radon measures on $\mathcal{O}$.
When $\mathcal{O}$ is a vector space, 
$\delta$ denotes the Dirac measure at the origin.

\subsection{The fundamental controlled solution}
\label{sec:fcs}

In this subsection we 
construct a ``fundamental controlled solution'' $v$
of the Schr{\"o}dinger equation on a segment controlled 
by Dirichlet conditions at both ends.
The precise definition is the following.

\begin{definition}
\label{def:fundcontrol}
The distribution $v\in C^{0}([0,T]; H^{-1}(]-L,L[))$
is a fundamental controlled solution
for the Schr{\"o}dinger equation
on $]0,T[\times ]-L,L[$ at cost $(A,\alpha)$
if 
\begin{align}
&i\partial_{t}v - \partial_{s}^{2}v  =  0
\quad {\rm in }\ \mathcal{D}'(]0,T[\times ]-L,L[)\ , 
\label{eqv1} \\
&v\res{t=0}  =  \delta \quad {\rm and }\quad v\res{t=T}  =  0 \ ,
\label{eqv2} \\
&\|v\|_{L^{2}(]0,T[; H^{-1}(]-L,L[))}
\leq Ae^{\alpha L^{2}/T } \ .
\label{eqv3}
\end{align}
\end{definition}

Theorem~\ref{th:1d} allows us to construct 
a family of fundamental controlled solutions depending on $L>0$ and $T>0$
with a good cost estimate 
thanks to the following proposition 
which shows that the upper bound for the controllability cost 
of the Schr{\"o}dinger equation on the segment $[0,L]$ 
controlled at one end
is the same as the controllability cost 
of the Schr{\"o}dinger equation on the twofold segment $[-L,L]$ 
controlled at both ends.

\begin{proposition}
\label{prop:twofold}
For any $\alpha > \alpha_{*}$ (cf. definition~\ref{defin:alpha}),
there exists $A>0$ such that, 
for all $L>0$, $T\in\, ]0,\inf(\pi/2,L)^{2}]$
and $v_{0}\in H^{-1}(-L,L)$,
there are $g_{-}$ and $g_{+}$ in $L^{2}(0,T)$
such that the solution $v\in C^{0}([0,\infty);H^{-1}(-L,L))$
of the following Schr{\"o}dinger equation on $\left[-L,L\right]$
controlled 
by $g_{-}$ and $g_{+}$:
\begin{equation} 
\label{eqHeattwofold}
i\partial_{t}v - \partial_{s}^{2} v=0
\quad {\rm in}\ ]0,T[\times ]-L,L[ ,\quad 
v\res{s=\pm L} =g_{\pm} ,\quad
v\res{t=0} = v_{0} 
\end{equation}
satisfies $v=0$ at $t=T$ and 
$\displaystyle
\|g_{\pm}\|_{L^{2}(0,T)} \leq Ae^{\alpha L^{2}/T }\|v_{0}\|_{H^{-1}(-L,L)}$.
\end{proposition}

\begin{proof}
By duality (cf.~\cite{DR77}), 
it is enough to prove the observation inequality:
$\displaystyle
\exists C>0, 
\forall v_{0}\in H^{1}_{0}(-L,L),
\|v_{0}\|_{H^{1}}
\leq C e^{\alpha L^{2}/T}
\|\partial_{s} e^{it\Delta}v_{0}{}_{\rceil s=\pm L}
\|_{L^{2}(0,T)^{2}}
$.
Applying theorem~\ref{th:1d} with $k=0$ to the odd part of $v_{0}$ 
and with $k=1$ to the even part of $v_{0}$
completes the proof 
(as in the proof of proposition~5.1 in~\cite{LMheatcost}).
\end{proof}

Applying proposition~\ref{prop:twofold} with 
$v_{0}=\delta\in H^{-1}(-L,L)$, 
and using Duhamel's formula to estimate $v$ in terms of $g_{\pm}=v\res{s=\pm L} $,
we obtain:

\begin{corollary}
\label{prop:fundcontrol}
For any $\alpha > \alpha_{*}$ (cf. definition~\ref{defin:alpha}),
there exists $A>0$ such that 
for all $L>0$ and $T\in\, ]0,\inf(\pi/2,L)^{2}]$
there is a fundamental controlled solution
for the Schr{\"o}dinger equation
on $]0,T[\times ]-L,L[$ at cost $(A,\alpha)$.
\end{corollary}

\subsection{The transmutation of waves controls into Schr{\"o}dinger controls}

In this subsection we perform a transmutation 
of a control for the wave equation 
into a control for the Schr{\"o}dinger equation.
Our transmutation formula (cf. (\ref{eqtrans})) 
can be regarded as the analogue of the formula (\ref{eqCGT})
with $F(\sigma)=\exp(it\sigma^{2})$
where the kernel $e^{-i\pi/4}e^{is^{2}/(4t)}/\sqrt{4\pi t}$,
which is the fundamental solution 
of the Schr{\"o}dinger equation on the line,
is replaced by the fundamental controlled solution 
that we have constructed in the previous subsection. 
To ensure existence of an exact control for the wave equation
we use the geodesics condition 
(cf. the footnotes on pages~\pageref{fnote:unique} and~\pageref{fnote:geod}):
\begin{theorem}[\cite{BLR92}]
\label{th:BLR}
Let $\Omega\Sset M$. 
Let $L_{\Omega}$ be the length of the longest generalized geodesic in $\adh{M}$
which does not intersect $\Omega$.
If $L>L_{\Omega}$ then
for all $(w_{0}, w_{1})$ and $(w_{2}, w_{3})$ 
in $H^{2}_{0}(M)\times H^{1}_{0}(M)$
there is a control function 
$f\in H^{1}_{0}(]0,L[\times M)$
such that the solution 
$w\in C^{0}([0,L];H^{2}_{0}(M))\cap C^{1}([0,L];H^{1}_{0}(M))$
of the mixed Dirichlet-Cauchy problem 
(n.b.~the time variable is denoted by $s$ here):
\begin{equation} 
\label{eqWave}
\partial_{s}^{2}w - \Delta w={\bf 1}_{\Omega_{L}} f 
\quad {\rm in}\ ]0,L[\times M, \quad 
w=0 \quad {\rm on}\ ]0,L[\times\partial M,
\end{equation}
with Cauchy data
$(w,\partial_{s}w)=(w_{0},w_{1})$ at $s=0$,
satisfies 
$(w,\partial_{s}w)=(w_{2},w_{3})$ at $s=L$.
Moreover, the operator 
$S_{W}: \left( H^{2}_{0}(M)\times H^{1}_{0}(M) \right)^{2}
\to H^{1}_{0}(]0,L[\times M)$
defined by $S_{W}\left((w_{0},w_{1}),(w_{2},w_{3})\right)=f$ 
is continuous.
\end{theorem}

From this theorem, 
the control transmutation deduces theorem~\ref{th:ub} for smooth data only:
\begin{proposition}
\label{prop:transmut}
For all $\alpha > \alpha_{*}$,  
there is exists $A>0$ such that for all 
$u_{0}\in H^{2}_{0}(M)$, $T\in]0,\min\{1,L_{\Omega}^{2}\}[$ 
and $L>L_{\Omega}$, 
there is a control $g\in L^{2}(\mathbb{R}\times M)$ which solves 
the controllability problem
\begin{gather}
  \label{equ1}
i\partial_{t}u - \Delta u=\cd{1}_{\Omega_{T}} g
\quad {\rm in }\  \mathcal{D}'(]0,T[\times M) 
 \quad {\rm and } \quad
u=0 \quad {\rm on}\ ]0,T[\times\partial M , \\
  \label{equ2}
u\res{t=0}  =  u_{0} \quad {\rm and }\quad u\res{t=T}  =  0 ,
\end{gather}
at cost 
$\|g\|_{L^{2}(\mathbb{R}\times M)}
\leq \|S_{W}\|  A e^{\alpha L^{2}/T }  \|u_{0}\|_{H^{2}_{0}(M)}$.
\end{proposition}
\begin{proof}
Let $\alpha > \alpha_{*}$, 
$T\in]0,\inf(1,L_{\Omega}^{2})[$ 
and $L>L_{\Omega}$ be fixed from now on.
Let $A>0$ and 
$v\in L^{2}(0,T;H^{-1}(]-L,L[))$ 
be the corresponding constant and fundamental controlled solution
given by corollary~\ref{prop:fundcontrol}.
We define $\vv\in L^{2}(\mathbb{R};H^{-1}(\mathbb{R}))$ 
as the extension of $v$ by zero,
i.e. $\vv(t,s)=v(t,s)$ on $]0,T[\times]-L,L[$
and $\vv$ is zero everywhere else.
It inherits from $v$ the following properties 
\begin{align}
&i\partial_{t}\vv - \partial_{s}^{2} \vv=0
\quad {\rm in }\  \mathcal{D}'(]0,T[\times ]-L,L[)\ , 
\label{eqvv1} \\
&\vv\res{t=0}  =  \delta \quad {\rm and }\quad \vv\res{t=T}  =  0 \ ,
\label{eqvv2} \\
&\|\vv\|_{L^{2}(\mathbb{R}; H^{-1}(\mathbb{R}))}
\leq Ae^{\alpha L^{2}/T } \ .
\label{eqvv3}
\end{align}

Let $u_{0}\in H^{2}_{0}(M)$ 
be an initial data for the Schr{\"o}dinger equation 
$(\ref{eqSchr})$.
Let $w$ and $f$
be the corresponding solution and control function 
for the wave equation obtained by applying theorem~\ref{th:BLR}
with $w_{0}=u_{0}$ and $w_{1}=w_{2}=w_{3}=0$.
Since theorem~\ref{th:BLR} still applies to 
any control time in $]L_{\Omega},L[$, 
we may assume that $f=0$ in a neighborhood of $s=0$. 
We define $\ww\in L^{2}(\mathbb{R};H^{2}_{0}(M))
\cap H^{1}(\mathbb{R};H^{1}_{0}(M))$
and $\ff\in H^{1}(\mathbb{R}\times M)$ 
as the extensions of $w$ and $f$ 
by reflection with respect to $s=0$,
i.e. $\ww(s,x)=w(s,x)=\ww(-s,x)$ 
and $\ff(s,x)=f(s,x)=\ff(-s,x)$ 
on $[0,L]\times M$, and by zero outside $[-L,L]\times M$.
Since $w_{1}=0$, equation $(\ref{eqWave})$ implies 
\begin{equation} 
\label{eqWaveww}
\partial_{s}^{2}\ww - \Delta \ww={\bf 1}_{]-L,L[\times \Omega} \ff 
\quad {\rm in}\  \mathcal{D}'(\mathbb{R}\times M), \quad 
\ww=0 \quad {\rm on}\ \mathbb{R}\times\partial M,
\end{equation}

The main idea of our proof is to use $\vv$
as a kernel to transmute $\ww$ and $\ff$
into a solution $u$ and a control $g$ 
for $(\ref{eqSchr})$. 
Since $\vv\in C^{0}(\mathbb{R};H^{-1}(\mathbb{R}))$, 
$\ww\in H^{1}(\mathbb{R};H^{1}_{0}(M))$ 
and $\ff\in H^{1}(\mathbb{R};L^{2}(M))$,
the transmutation formulas 
\begin{equation} 
\label{eqtrans}
u(t,x)=\int_{\mathbb{R}}\vv(t,s)\ww(s,x)\, ds  
 \quad {\rm and } \quad 
g(t,x)=\int_{\mathbb{R}} 
\vv(t,s)\ff(s,x)\, ds 
\ ,
\end{equation}
define functions 
$u\in C^{0}(\mathbb{R};H^{1}_{0}(M))$ 
and $g\in L^{2}(\mathbb{R}\times M)$.
Since $\ww(s,x)=\partial_{s}\ww(s,x)=0$ for $|s|=L$,
equations $(\ref{eqWaveww})$ and 
$(\ref{eqvv1})$ imply $(\ref{equ1})$.
The property $(\ref{eqvv2})$ of $\vv$ implies $(\ref{equ2})$.
Since 
$\|g\|_{L^{2}(\mathbb{R}\times M)}
\leq \|\vv\|_{L^{2}(\mathbb{R}; H^{-1}(\mathbb{R}))}
\|\ff\|_{H^{1}(\mathbb{R};L^{2}(M))}$,
the estimates $(\ref{eqvv3})$ and
$\|\ff\|_{H^{1}(\mathbb{R}\times M)}
\leq \sqrt{2}\|S_{W}\|\, \|u_{0}\|_{H^{2}_{0}(M)}
$
complete the proof of proposition~\ref{prop:transmut}.
\end{proof}

\begin{proof}[Proof of theorem~\ref{th:ub}]
Let $\alpha > \alpha_{*}$, $L>L_{\Omega}$ and $\eps\in ]0,1[$.

According to corollary~\ref{cor:smoothing}, there are $d>0$ and $C>0$
such that for all $T\in]0,\inf(1,L_{\Omega}^{2})[$,
there is a $T'\in [\eps T/2,\eps T]$ such that the following 
smoothing property at low control cost holds: 
for any initial data $u_{0}\in L^{2}(M)$ 
there is a $g_{1}\in L^{2}(\Omega_{T'})$ such that 
$\|g_{1}\|_{L^{2}}\leq \|u_{0}\|_{L^{2}}C/\sqrt{T'}$
and the solution $u$ of $(\ref{eqSchr})$ with control $g_{1}$ satisfies
$\|u\res{t=T'}\|_{H^{2}}
\leq C(T')^{-2}\|u_{0}\|_{L^{2}}$.

According to proposition~\ref{prop:transmut}, 
there is a $g_{2}\in L^{2}(]T',T[\times\Omega)$ such that 
$\|g_{2}\|_{L^{2}}\leq 
A \|S_{W}\|  e^{\alpha L^{2}/(T-T')}  \|u\res{t=T'}\|_{H^{2}}$
and the solution $u$ of $(\ref{eqSchr})$ 
with the control $g$ obtained by applying successively 
$g_{1}$ and $g_{2}$ satisfies $u\res{t=T'}=0$.

These estimates with $T'\in [\eps T/2,\eps T]$  
combine into a cost estimate for $g$:
\begin{gather*}
\|g_{2}\|_{L^{2}}^{2}=\|g_{1}\|_{L^{2}}^{2}+\|g_{2}\|_{L^{2}}^{2}
\leq  \|u_{0}\|_{L^{2}}^{2}\left(
\frac{2C^{2}}{\eps T} + \frac{4C}{\eps^{2} T^{2}}A^{2} \|S_{W}\|^{2}
e^{2\alpha L^{2}/(1-\eps)T}
\right)
\ . \end{gather*}

With the dual definition of $C_{T,\Omega} $ 
given after definition~\ref{defin:cost}, this estimate proves: 
$\displaystyle
\limsup_{T\to 0} T \ln C_{T,\Omega} \leq \alpha L^{2}/(1-\eps)$.
Letting $\alpha$, $L$ and $\eps$  tend respectively to 
$\alpha_{*}$, $L_{\Omega}$ and $0$ 
completes the proof of $(\ref{eq:ub})$.
\end{proof}

\section{Upper bound for some examples violating the geodesics condition}
\label{sec:ubprod}

In this section, 
we deduce from theorem~\ref{th:1d} and~\ref{th:ub} 
that the same upper bounds are satisfied
for some Schr{\"o}dinger evolution groups of product type
violating the geodesics condition.
The proof elaborates on the yet unpublished remark of Burq 
(back in 1992, cf~\cite{BZbb}) 
that the result of~\cite{Har89} can be extended to product manifolds
with a much simpler proof:
the point here is that the controllability cost is tracked.

The following lemma generalizes this remark to the abstract setting 
for the theory of observation and control (cf. \cite{DR77}). 

\begin{lemma}
\label{lem:abstprod}
Let $X$, $Y$ and $Z$ be Hilbert spaces
and $I$ denote the identity operator on each of them. 
Let $A:D(A)\to X$ and $B:D(B)\to Y$ be generators of
strongly continuous semigroups of bounded operators on $X$ and $Y$.
Let $C:D(C)\to Z$ be a densely defined operator on $X$
such that $e^{tA}D(C)\subset D(C)$ for all $t>0$.
Let $X\ctensprod Y$ and $Z\ctensprod Y$ denote the closure of 
the algebraic tensor products $X\otensprod Y$ and $Z\otensprod Y$ 
for the natural Hilbert norms.
The operator $C\otensprod I:D(C)\otensprod Y \to Z\ctensprod Y$ 
is densely defined on $X\ctensprod Y$.

i) The operator $A\otensprod I+I\otensprod B$ defined on 
the algebraic $D(A)\otensprod D(B)$ is closable
and its closure, denoted $A+B$, 
generates a strongly continuous semigroup of bounded operators on $X\ctensprod Y$
satisfying:
\begin{gather}
\label{eqnormprod}
\forall t\geq 0, \forall (x,y)\in D(C)\times Y,\quad 
\norm{(C\otensprod I) e^{t(A+B)}(x\otensprod y)}=\norm{Ce^{tA}x}\, \norm{e^{tB}y}
\end{gather}

ii) If $iB$ is self-adjoint, then for all $T\geq 0$: 
\begin{gather}
\label{eqcostprod}
\inf_{\psi\in X\ctensprod Y, \norm{\psi}=1} \int_{0}^{T}\norm{(C\otensprod I)e^{t(A+B)}\psi}^{2} dt
= \inf_{x\in X, \norm{x}=1} \int_{0}^{T}\norm{Ce^{tA}x}^{2} dt
\ .
\end{gather}
\end{lemma}

\begin{proof}
Let $G$ denote the generator of 
the strongly continuous semigroup $t\mapsto e^{tA}\otensprod e^{tB}$
(defined since the natural Hilbert norm is a uniform cross norm, cf. \cite{Sch50}).
Since $D(A)\otensprod D(B)$ is dense in $X\otensprod Y$ and invariant by $t\mapsto e^{tG}$,
it is a core for $G$ (cf. theorem~X.49 in~\cite{RS}).
Since $A\otensprod I+I\otensprod B=G\res{D(A)\otensprod D(B)}$, 
it is closable and $A+B=G$.
Therefore $e^{t(A+B)}=e^{tA}\otensprod  e^{tB}$ and (\ref{eqnormprod}) follows 
(by the cross norm property).

To prove point~ii), we denote the left and right hand sides of (\ref{eqcostprod})
by $\mathcal{I}_{A+B}$ and $\mathcal{I}_{A}$.
Taking $\psi=x\otensprod  y$ with $\norm{y}=1$, 
$\mathcal{I}_{A+B}\leq\mathcal{I}_{A}$ results from (\ref{eqnormprod}).
To prove $\mathcal{I}_{A+B}\geq\mathcal{I}_{A}$, 
we only consider the case in which both $X$ and $Y$ are 
infinite dimensional and separable
(this simplifies the notation and the other cases are similar).
Let $(e_{n})_{n\in \N}$ and $(f_{n})_{n\in \N}$ be orthonormal bases for $X$ and $Y$.
Since $(e_{n}\otensprod  f_{m})_{n,m\in \N}$ is an orthonormal base for $X\ctensprod Y$,
any $\psi\in X\ctensprod Y$ writes: 
\begin{gather*}
\psi=\sum_{m}x_{m}\otensprod  f_{m} 
\quad \text{with }\ x_{m}=\sum_{n}c_{n,m}e_{n}
\ \text{ and }\ 
\norm{\psi}^{2}=\sum_{n,m}\abs{c_{n,m}}^{2}=\sum_{m}\norm{x_{m}}^{2}
\ .
\end{gather*}
Since $iB$ is self-adjoint, 
$t\mapsto e^{tB}$ is unitary for all $t\geq 0$
so that  $(e^{tB}f_{n})_{n\in \N}$ is orthonormal.
Therefore, using (\ref{eqnormprod}):
\begin{gather*}
\norm{Ce^{t(A+B)}\psi}^{2}=\norm{\sum_{m}(Ce^{tA}x_{m})\otensprod (e^{tB}f_{m})}^{2}
=\sum_{m}\norm{Ce^{tA}x_{m}}^{2}
\ .
\end{gather*}
By definition, 
$\int_{0}^{T}\norm{Ce^{tA}x_{m}}^{2} dt\geq \mathcal{I}_{A}\norm{x_{m}}^{2}$.
Summing up over $m\in\N$, we obtain:
\begin{gather*}
\int_{0}^{T}\norm{(C\otensprod  I)e^{t(A+B)}\psi}^{2} dt
=\int_{0}^{T}\sum_{m}\norm{Ce^{tA}x_{m}}^{2}
\geq \mathcal{I}_{A}\sum_{m}\norm{x_{m}}^{2}
= \mathcal{I}_{A}\norm{\psi}^{2}
\ .
\end{gather*}
This proves $\mathcal{I}_{A+B}\geq\mathcal{I}_{A}$ 
and completes the proof of lemma~\ref{lem:abstprod}.
\end{proof}

\begin{remarks}
When $C$ is an admissible observation operator, 
(\ref{eqcostprod}) says that 
the cost of observing $t\mapsto e^{t(A+B)}$ through $C\otensprod  I$ in time $T$
is exactly 
the cost of observing $t\mapsto e^{tA}$ through $C$ in time $T$.

If $A$ and $B$ are self-adjoint, then $A+B$ defined in lemma~\ref{lem:abstprod}
is self-adjoint (cf. theorem~VIII.33 in~\cite{RS}).

The proof of part~i) of lemma~\ref{lem:abstprod} is still valid 
if $X$, $Y$ and $Z$ are Banach spaces
and $X\ctensprod Y$ and $Z\ctensprod Y$ are closures 
with respect to some uniform cross norms (cf. \cite{Sch50}).
\end{remarks}

Theorem~\ref{th:ubprod} is a particular case of the following
direct consequence of lemma~\ref{lem:abstprod} and theorem~\ref{th:1d}
(with $X=Z=L^{2}(M)$, $Y=\mathcal{B}$, $A=i\Delta$ and bounded $C=\cd{1}_{\Omega}$).

\begin{theorem}
\label{th:prod}
Let $B$ be a self-adjoint operator on a Hilbert space $\mathcal{B}$.
The operator $H=\Delta\otimes\id_{\mathcal{B}}+\id_{L^{2}(M)}\otimes B$ 
is essentially self-adjoint on $\mathcal{H}=L^{2}(M)\otimes\mathcal{B}$.
For all $T>0$ and $\Omega\Sset M$, 
$C_{T,\Omega}$ (cf. definition~\ref{defin:cost}) 
is also the cost of controlling the Schr{\"o}dinger group 
$t\mapsto e^{itH}$ on $\mathcal{H}$
with controls in $L^{2}(\Omega)\otimes\mathcal{B}$, i.e. 
$C_{T,\Omega}$ is the best constant in the observability inequality:
$\forall v\in \mathcal{H}$,
$\norm{v}_{\mathcal{H}}
\leq C_{T,\Omega}\norm{\cd{1}_{\Omega}e^{itH}v}_{L^{2}(]0,T[;\mathcal{H})}
$.
In particular: 
\begin{equation*}
\limsup_{T\to 0} T\ln 
\sup_{ v\in \mathcal{H}\setminus\set{0} }
\frac{\|v\|_{\mathcal{H}}}
{\norm{\cd{1}_{\Omega}e^{itH}v}_{L^{2}(]0,T[;\mathcal{H})}}
\leq \alpha_{*}L_{\Omega}^{2}
\ , 
\end{equation*}
where 
$\alpha_{*}$ is defined in (\ref{eqalphastar})
and $L_{\Omega}$ is the length of the longest generalized geodesic in $\adh{M}$
which does not intersect $\Omega$.
\end{theorem}

\begin{remarks}
Note that $\mathcal{H}=L^{2}(\Omega;\mathcal{B})$ and 
$L^{2}(0,T;\mathcal{H})=L^{2}(]0,T[\times M;\mathcal{B})$
when $\mathcal{B}$ is separable (cf. theorem~II.10 in~\cite{RS}).
With $\mathcal{B}=L^{2}(\Mt)$, $B=\Deltat$, $\mathcal{H}=L^{2}(M\times\Mt)$,
theorem~\ref{th:prod} proves theorem~\ref{th:ubprod}.

The semi-internal controllability of a rectangular plate proved in~\cite{Har89}
corresponds to the setting $M=[0,X]$, 
$\mathcal{B}=L^{2}([0,Y])$, $B=\partial_{y}^{2}$ with Dirichlet condition, 
$\mathcal{H}=L^{2}([0,X]\times [0,Y])$.
Note that our theorem still applies 
to an infinite strip $[0,X]\times \R$
controlled from any infinite strip  $[a,b]\times \R$
with $[a,b]\subset ]0,X[$.

The resolvent method introduced in~\cite{BZbb} also yields 
the controllability in theorem~\ref{th:prod} for some control time
(and for any positive control time by a temporal black box), 
but it does not keep track of the cost.
\end{remarks}

The following analogue of theorem~\ref{th:prod}
for the boundary controllability of cylinders from one end
is a direct consequence of lemma~\ref{lem:abstprod} and theorem~\ref{th:1d} 
(with $X=H^{k}_{A}(0,X)$, $Y=\mathcal{B}$, $Z=\R$, $D(C)=D(A)$ 
and $Cu=\partial_{x}u\res{x=X}$).

\begin{theorem}
\label{th:cyl}
Let $B$ be a self-adjoint operator on a Hilbert space $\mathcal{B}$.
Let $A$ be the Sturm-Liouville operator on $L^{2}(0,X)$
and $L$ be the length of $[0,X]$ 
defined in subsection \ref{sec:1d}.
The operator $H=A\otimes\id_{\mathcal{B}}+\id_{L^{2}(0,X)}\otimes B$ 
is essentially self-adjoint on $\mathcal{H}=L^{2}(0,X)\otimes\mathcal{B}$.
For any $\alpha > \alpha_{*}$ defined by $(\ref{eqalphastar})$,
there exists $C>0$ such that, 
for any coefficients $(\ref{eqSLAssum})$, 
setting $k=1$ if $b_{1}=0$ and $k=0$ otherwise,
for all $T\in\, ]0,\inf(\pi,L)^{2}]$:
\begin{gather*}
\forall v\in \mathcal{H}^{k}=H^{k}_{A}(0,X)\otimes\mathcal{B},\quad
\|v\|_{\mathcal{H}^{k}} 
\leq  C\exp(\alpha L^{2}/T)
\|\partial_{x}^{k}e^{itH}v\res{x=X}\|_{L^{2}(0,T;\mathcal{B})}
\ .
\end{gather*}
\end{theorem}

\begin{remarks}
With $\mathcal{B}=L^{2}(\Mt)$, $B=\Deltat$, $\mathcal{H}=L^{2}(C)$,
$A=\partial_{x}^{2}$, $k=1$,
this theorem applies to 
the Schr{\"o}dinger equation on the cylinder
$C=[0,X]\times \Mt$ controlled at the end $\Gamma=\{X\}\times \Mt$,
with a base $\Mt$ as in theorem~\ref{th:prod}.
The segment $S=[0,X]$ is endowed with a Riemannian metric,
$L$ denotes the total length of $[0,X]$ 
and $\Delta_{S}$ denotes the Dirichlet Laplacian on $[0,X]$,
so that the Laplacian on the $(n+1)$-dimensional product manifold $C$ is
$\Delta_{C}=\Delta_{S}+\Deltat$.
In this setting, the controllability cost is the best constant, 
denoted $C_{T,\Gamma}$, in the observation inequality:
\begin{gather}
\label{eq:cylcost}
\forall u_{0}\in \mathcal{H}^{1}=H^{1}_{0}(S;L^{2}(\Mt)),\quad
\|u_{0}\|_{\mathcal{H}^{1}} 
\leq C_{T,\Gamma}
\|\partial_{s}e^{it\Delta_{C}}u_{0}{}\res{\Gamma}\|_{L^{2}(]0,T[\times M)} \ .
\end{gather}
Although the geodesics condition is not satisfies for $\Gamma$ in $C$,
theorem~\ref{th:cyl} proves 
that the controllability cost $C_{T,\Gamma}$ satisfies,
as in theorem~\ref{th:1d},
an upper bound of the same type as 
the lower bound in theorem~\ref{th:lb}:
$\forall \alpha>\alpha_{*}$, $\exists \beta>0$, 
$C_{T,\Gamma}\leq \beta\exp(\alpha L^{2}/T)$.

Note that the observability inequality (\ref{eq:cylcost}) 
does not hold in the ``~energy space~'', i.e. 
the space $\mathcal{H}^{1}=H^{1}_{0}(S,L^{2}(\Mt))$ 
cannot be replaced by $H^{1}_{0}(C)$.

The boundary controllability of a rectangular plate from one side   
was proved in~\cite{KLS85} (theorem~2).
\end{remarks}



\providecommand{\bysame}{\leavevmode\hbox to3em{\hrulefill}\thinspace}
\providecommand{\MR}{\relax\ifhmode\unskip\space\fi MR }
\providecommand{\MRhref}[2]{%
  \href{http://www.ams.org/mathscinet-getitem?mr=#1}{#2}
}
\providecommand{\href}[2]{#2}

\end{document}